\documentclass[12pt]{interact}

\usepackage{graphics,xcolor,graphicx,subfigure,amssymb,adjustbox,a4wide}

\newtheorem{thm}{Theorem}[section]
\newtheorem{remark}[thm]{Remark}
\newtheorem{lemma}[thm]{Lemma}
\newtheorem{defi}{Definition}[section]
\newtheorem{conjecture}[thm]{Conjecture}
\def\xsum{\mathop{\sum\nolimits'}}

\begin{document}

\title{Interplay between critical and off-critical zeros of two-dimensional Epstein zeta functions}


\author{
\name{Laurent B\'{e}termin}
\affil{Institut Camille Jordan, Universit\'e Claude Bernard Lyon 1,
69622 Villeurbanne, France} \medskip
\name{Ladislav \v{S}amaj \and Igor Trav\v{e}nec} 
\affil{Institute of Physics, Slovak Academy of Sciences, 
D\'ubravsk\'a cesta 9, 84511 Bratislava, Slovakia}}

\date{Received:   / Accepted: }

\maketitle

\begin{abstract}
The two-dimensional Epstein zeta function associated to a rectangular
lattice with spacings $a_x=1$ and $a_y=\Delta\in \mathbb{R}_+^*$, defined by
$\zeta^{(2)}(s,\Delta) = \frac{1}{2} \xsum_{j,k} (j^2+\Delta^2 k^2)^{-s}$
$(\Re(s)>1)$, where the sum goes over all integers except of the origin
$(j,k)=(0,0)$, is studied.
It can be analytically continued to the whole complex $s$-plane
except for the point $s=1$.
The nontrivial zeros $\{ \rho=\rho_x+i\rho_y \}$ of the Epstein zeta
function, defined by $\zeta^{(2)}(\rho,\Delta)=0$, split into ``critical''
zeros (on the critical line $\rho_x=\frac{1}{2}$) and ``off-critical'' zeros
($\rho_x\ne\frac{1}{2}$).
This work presents rigorous asymptotic and analytic results as well as
numerical investigations.
According to the present numerical calculation, the critical zeros form
open or closed curves $\rho_y(\Delta)$ in the plane $(\Delta,\rho_y)$.
Nearest critical zeros merge at special points
$(\Delta^*,\rho_y^*)$, referred to as left/right edge zeros,
which are defined by a divergent tangent ${\rm d}\rho_y/{\rm d}\Delta$.
Each of these critical edge zeros gives rise to a continuous curve of
off-critical zeros which can thus be generated systematically.
As a rule, each curve of off-critical zeros joins a pair of left and right
edge zeros.
It is shown that in the regions of small/large values of the anisotropy
parameter $\Delta$ the Epstein zeta function can be approximated adequately
by a function which reveals an equidistant distribution of critical zeros
along the imaginary axis in the limits $\Delta\to 0$ and $\Delta\to\infty$.
It is also numerically found that for each
$\Delta\in (0,\Delta_c^*]\cup [1/\Delta_c^*,\infty)$ with
$\Delta_c^*\approx 0.141733$ there exists a pair of \emph{real} off-critical
zeros, their $\rho_x$ components going to the borders $0$ and $1$ of
the critical region in the limits $\Delta\to 0,\infty$.

\end{abstract}
  
  
\begin{amscode}
11E45, 11M41, 11R42
\end{amscode}

\renewcommand{\theequation}{1.\arabic{equation}}
\setcounter{equation}{0}

\section{Introduction} \label{Sec1}

\subsection{The zeta functions and their zeros.} We consider the two-dimensional (2D) Epstein zeta function \cite{Epstein03,Epstein07} associated to the rectangular lattice $a_x\mathbb{Z}\oplus a_y\mathbb{Z}$ with spacings $(a_x,a_y)=(1,\Delta)$, $\Delta\in \mathbb{R}_+^*$, defined by
\begin{equation} \label{sum1}
\zeta^{(2)}(s,\Delta) := \frac{1}{2} \xsum_{(j,k)\in \mathbb{Z}^2}
\frac{1}{(j^2+\Delta^2 k^2)^s} \qquad \Re(s)>1,
\end{equation}
where $\xsum$ means that the term $(j,k)=(0,0)$ is excluded from the summation. This represents a natural 2D extension of the one-dimensional (1D) Riemann zeta function \cite{Riemann1859} defined by
\begin{equation} \label{RiemannZeta}
\zeta(s)= \frac{1}{2} \xsum_{j\in \mathbb{Z}}\frac{1}{\vert j\vert^{s}}
= \sum_{j\in \mathbb{N}} \frac{1}{j^{s}}.
\end{equation}
The sum \eqref{sum1} naturally arises in Physics as the energy per point of a system of identical particles located on the lattice sites of $\mathbb{Z}\oplus \Delta\mathbb{Z}$ and interacting pairwisely through the Riesz potential $r\mapsto 1/r^{2s}$ \cite{Brau11}. Notice that the prefactor $\frac{1}{2}$ is present because each interaction
energy is shared by a pair of particles.

The function $\zeta^{(2)}(s,\Delta)$ possesses the obvious symmetry
\begin{equation} \label{Deltasymmetry}
\zeta^{(2)}(s,\Delta) = \frac{1}{\Delta^{2s}} \zeta^{(2)}(s,1/\Delta) 
\end{equation}
which means that the values of $\Delta$ can be constrained to either of
the intervals $(0,1]$ or $[1,\infty)$. Similarly as the Riemann zeta function,
the Epstein zeta function can be analytically continued to the whole complex
$s$-plane \cite{Elizalde89a,Ennola64,Epstein03,Epstein07}.
The region $0<\Re(s)<1$ is referred to as the critical strip, the critical
line is defined by $\Re(s)=\frac{1}{2}$.

For $\Delta^2\in \{1,2,3,4,7\}$, the 2D lattice sum (\ref{sum1})
can be expressed as a product of simpler 1D sums,
namely Dirichlet $L$-series, whereas for other special
integer values of $\Delta$ it is expressible as \emph{a sum} of products
of Dirichlet $L$-functions \cite{Borwein13}.
For instance, in the isotropic square lattice case $\Delta^2=1$, it holds
\begin{equation} \label{zeta1}
\zeta^{(2)}(s,1) = 2 \zeta(s) \beta(s), \quad \textnormal{where}\quad \beta(s)
:= \sum_{j=0}^{\infty} \frac{(-1)^j}{(2j+1)^s} =
\frac{1}{2^{2s}} \left[ \zeta\left( s,\frac{1}{4}\right)  -
\zeta\left( s,\frac{3}{4}\right) \right]
\end{equation}  
is the Dirichlet beta function and
$\zeta(s,a):= \sum_{j=0}^{\infty} (j+a)^{-s}$ denotes
the Hurwitz zeta function.

The zeros $\{\rho=\rho_x+i\rho_y\}$ of the Epstein zeta function
$s\mapsto\zeta^{(2)}(s,\Delta)$ associated to the parameter $\Delta$ are defined by the equality $\zeta^{(2)}(\rho,\Delta)=0$.
Besides the trivial zeros at $\rho\in \mathbb{Z}_-$ there exist two kinds
of nontrivial zeros: 
\begin{itemize}
\item the ``critical'' zeros which lie on the critical line
$\rho_x=\frac{1}{2}$;
\item the ``off-critical'' zeros which have
$\rho_x\ne \frac{1}{2}$.
\end{itemize}
The Riemann hypothesis states that all nontrivial zeros of the 1D Riemann zeta
function $\zeta$ are constrained to the corresponding critical line
$\rho_x=\frac{1}{2}$ \cite{Riemann1859}.
Provided that the Riemann hypothesis holds, all nontrivial zeros of
the Dirichlet beta function (\ref{zeta1}) are constrained to
the critical line $\rho_x=\frac{1}{2}$ as well \cite{Lander18}.
Consequently, all nontrivial zeros of the Epstein zeta function
$\zeta^{(2)}(s,1)$ associated to the square lattice,
given by (\ref{zeta1}), lie on the critical line.
Similar phenomenon is expected also for $\Delta^2\in\{2,3,4,7\}$ when the
Epstein zeta function factorizes itself into the product of Dirichlet
$L$-functions; for analytic and numerical studies of the statistics of gaps
between critical zeros, see \cite{Baier17,Bogomolny94,Bombieri87,Jutila05}.
This is no longer true for anisotropic (i.e. non-square) rectangular
lattices (\ref{sum1}) with other integer values of $\Delta^2$.
The first off-critical zero of $\zeta^{(2)}(s,\Delta)$ was detected for
$\Delta^2=5$ \cite{Potter35}, since then many other off-critical zeros were
identified \cite{Bateman64,Davenport36,McPhedran16,Stark67}.

\subsection{Main results and related works}
\begin{figure}[t]
\centering
\includegraphics[clip,width=0.75\textwidth]{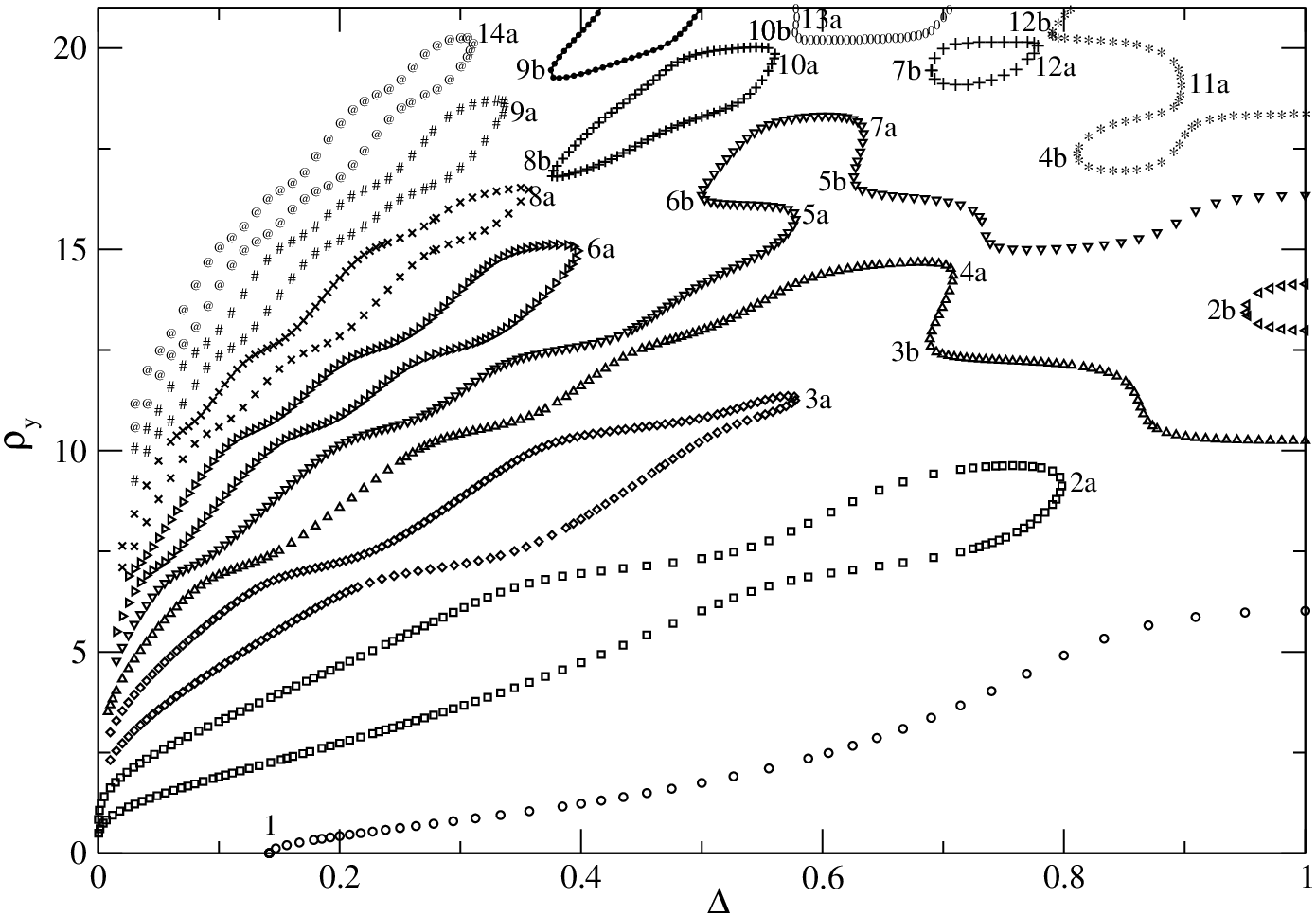}
\caption{Numerical results for the critical zeros in the region of
$0<\Delta\le 1$ and $0\le \rho_y\le 21$.
Zeros lying on the same (closed or open) curve are denoted by
a common open symbol (circle, square, triangle,$\ldots$).  
The right and left edge points are denoted as 
2a, 3a, 4a,$\ldots$ and 1, 2b, 3b, 4b,$\ldots$, respectively.}
\label{allcritintro}
\end{figure}

\begin{figure}[t]
\centering
\includegraphics[clip,width=0.75\textwidth]{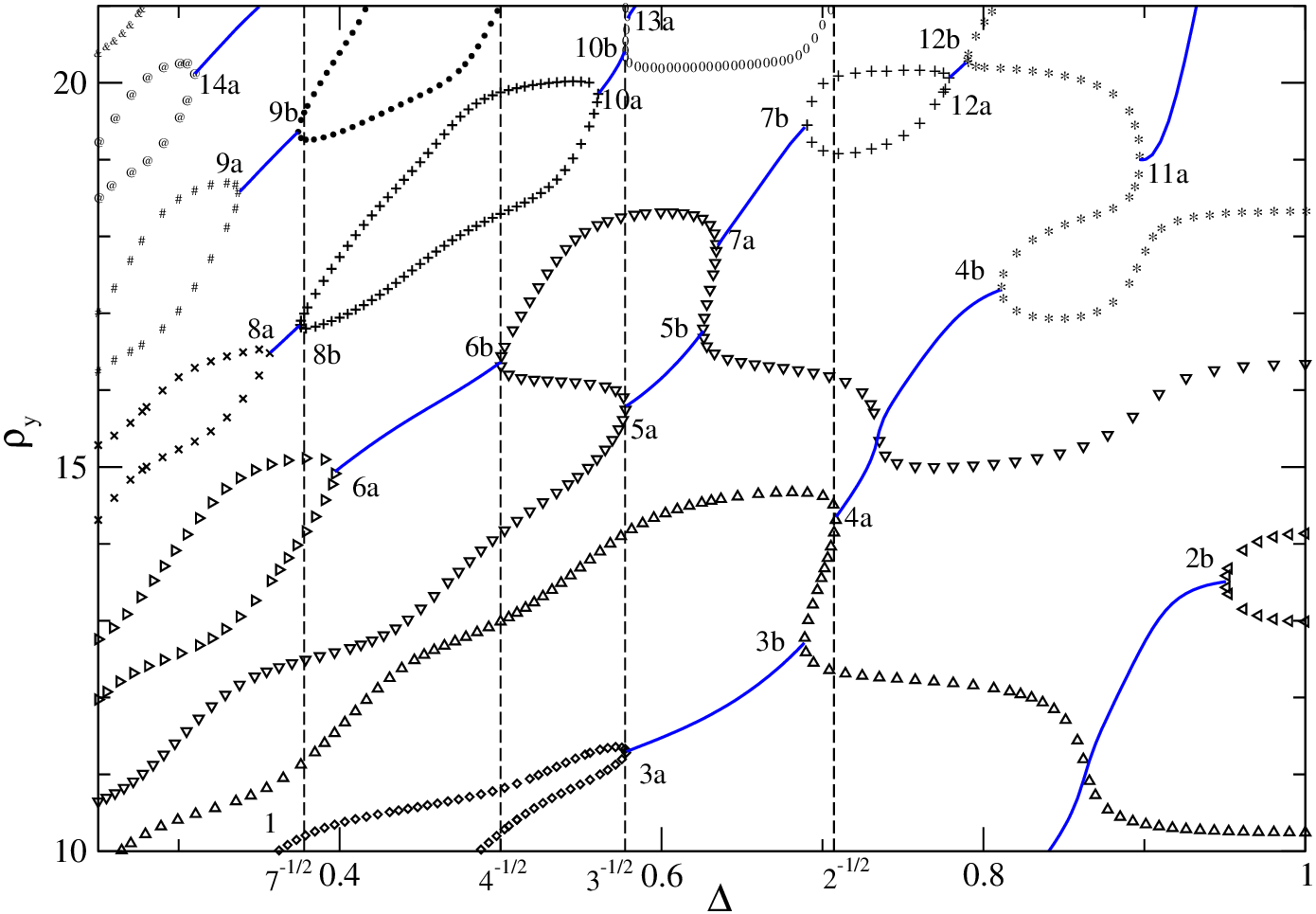}
\caption{
A fragment of Figure \ref{allcritintro}.  
Numerical results for solid curves of off-critical zeros interconnecting
the critical right (notation ``a'') and left (``b'') edge zeros.
The evolution of the component $\rho_x$ along the solid curves of off-critical
zeros is not indicated.
The vertical dashed lines correspond to the values of
$\Delta^2=1,2,3,4,7$ with no off-critical zeros.}
\label{critoffintro}
\end{figure}

In this paper, we perform a numerical investigation of the critical zeros of
$s\mapsto \zeta^{(2)}(s,\Delta)$ for $\Delta \in (0,1]$ in
a compact set and explain how the off-critical zeros of
this Epstein zeta function emerge from certain critical ``edge" zeros.
Most of our numerical findings are based on rigorous asymptotic and analytic
results.
Furthermore, we systematically check the degree of accuracy between
our approximations and the direct numerical computations from
more complicated equations.
More precisely, our observation and results are the following:
\begin{enumerate}
\item For $\zeta^{(2)}(s,\Delta)$ one obtains numerically in the plane $(\Delta,\rho_y)$ continuous curves of critical zeros $\rho_y(\Delta)$ (remember that $\rho_x=1/2$) and identifies critical edge zeros on them that correspond to the points of the curve with vertical tangents. These curves and critical edge zeros are depicted in Figure \ref{allcritintro} and their construction is explained in Section \ref{Sec3}.
\item The critical edge zeros are the starting points to
generate systematically (here numerically) continuous curves of off-critical
zeros. These continuous curves, connecting a pair of left and
right edge zeros, are depicted in Figure \ref{critoffintro} and their
construction is explained in Section \ref{Sec4}.
\item The analysis of the limits $\Delta\to 0$ and $\Delta\to \infty$ reveals
an equidistant distribution of critical zeros along the imaginary axis,
with spacing between the nearest zeros going to zero as
$\mathcal{O}(1/\vert\log \Delta\vert)$.
This analysis is done in Section \ref{Sec33}.
\item A pair of real off-critical zeros is numerically found for each
$\Delta\in (0,\Delta_c^*]\cup [1/\Delta_c^*,\infty)$ with
$\Delta_c^*\approx 0.141733$.
This is explained in Section \ref{Sec42}. Furthermore, our asymptotic study
for small values of $\Delta$ combined with a result by Montgomery \cite{Mont}
on theta functions (see also \cite{FaulhuberSteinerberger17}) yields
the following conjecture that we have also numerically checked.
\end{enumerate}
\begin{conjecture}\label{conjecture}
The unique solution of $\zeta^{(2)}(1/2,\Delta)=0$ for
$\Delta\in (0,1)$ is $\displaystyle \Delta=\Delta_c^* = \frac{e^\gamma}{4\pi}$
where $\gamma$ is the Euler-Mascheroni constant.
\end{conjecture}

Our work is motivated by another recent paper on $d$-dimensional Epstein zeta
function where surprising similarities appear.
Indeed, the Epstein zeta function can be defined for any lattice structure
and the two last authors have recently studied in \cite{Travenec22}
the properties of the Epstein zeta function on an isotropic hypercubic lattice
\begin{equation} \label{Epsteind}
\zeta^{(d)}(s) = \frac{1}{2} \xsum_{(n_1,n_2,\ldots,n_d)\in \mathbb{Z}^d}
\frac{1}{(n_1^2+n_2^2+\cdots+n_d^2)^{s/2}} \qquad \Re(s)>d , 
\end{equation}
where $d=1,2,\ldots$ is the spatial dimension and the Riesz interaction
energy between two particles at distance $r$ was chosen as $1/r^s$
to reproduce in $d=1$ the Riemann zeta function, $\zeta^{(1)}(s)=\zeta(s)$.
An analytic continuation of the Epstein zeta function $\zeta^{(d)}(s)$
to the whole complex $s$-plane was constructed; the corresponding formula
is applicable for the spatial dimension $d$ being a continuous variable
ranging from $0$ to $\infty$.
Numerical calculations of critical zeros (with $\rho_x=\frac{d}{2}$) indicate
that they form closed or semi-open curves in the plane
$\left( \rho_x=\frac{d}{2},\rho_y\right)$.
Each curve involves a number of left/right ``edge'' points
at which a couple of nearest critical zeros merge.
The critical edge zeros give rise to two tails of off-critical
zeros, coupled via a symmetry, with continuously varying dimension $d$
and the $\rho_x$-component along each tail.
This fact permits one to avoid a ``blind'' search for off-critical zeros,
the off-critical zeros are generated systematically starting
from the critical edge zeros.
Another benefit of the method is an exact treatment of the limits
$d\to 0$ and $d\to \infty$.
An exact formula for $\lim_{d\to 0} \zeta^{(d)}(s)/d$ was derived.
In the limit $d\to\infty$, an equidistant distribution of critical zeros
along the imaginary axis was obtained, with spacing between the
nearest zeros going to zero as $2\pi/\log(d)$.
As a by-product of the formalism, a conjugate pair of \emph{real}
off-critical zeros was found for each dimension $d>9.24555\ldots$.

The problem of determining zeros for the $d$-dimensional isotropic Epstein
zeta function (\ref{Epsteind}) seems at first sight to be unrelated to that
for the 2D Epstein zeta function (\ref{sum1}).
In analogy with \cite{Travenec22}, each of the curves of critical zeros
in the $(\Delta,\rho_y)$ plane pictured in Figure
\ref{allcritintro} contains a number of left/right edge zeros.
The mechanism of generation of off-critical zeros from these edge points,
depicted in Figure \ref{critoffintro} and explained in Section \ref{Sec41},
is the same as in the previous case of the $d$-dimensional isotropic
Epstein zeta function.

However, there are small differences in the form of curves of off-critical
zeros.
In the present 2D case, there is only one curve going out of
the classical edge zero, at each point of the curve there are two conjugate
off-critical zeros with $\rho_x$ and $1-\rho_x$ coordinates and each
curve of off-critical zeros connects just one left and one right edge points.
In the previous $d$-dimensional case, there are two tails of
off-critical zeros starting at each edge point, coupled via a symmetry, 
and a tail either interconnects two edge zeros or starts from one edge zero and
ends at some unspecial point in $d=0$ or $d\to\infty$ dimensions.

Another analogy with the previous paper \cite{Travenec22} is an equidistant
distribution of critical zeros along the imaginary axis in special regions
of model's parameters.
In the case of the $d$-dimensional isotropic Epstein function (\ref{Epsteind}),
the critical zeros are distributed equidistantly in the limit $d\to\infty$
with the spacing $2\pi/\log(d)$ between the nearest-neighbor zeros.
In the present case of the 2D anisotropic Epstein function (\ref{sum1}),
the critical zeros are distributed equidistantly in the limits $\Delta\to 0$
and $\Delta\to\infty$, with the spacing of order $\pi/\vert\log(\Delta)\vert$
between the nearest-neighbor zeros (see Section \ref{Sec33}).

The next similarity with \cite{Travenec22} consists in the appearance of
a couple of real off-critical zeros with the component $\rho_y=0$,
discussed in Section \ref{Sec42}.

We would like to emphasize that the presented mechanism of generation
of continuous curves of off-critical zeros from the edge critical ones   
might be not the only possible one.
Our crucial assumption was that the deviation of $\rho_x$ from its critical
value $\frac{1}{2}$ changes continuously when the off-critical zero
goes out the edge zero.
A discontinuous change of $\rho_x$ was excluded from our analysis. 

\medskip

\textbf{Plan of the paper.}
Section \ref{Sec2} concerns technicalities, like recalling the analytic
continuation of the 2D Epstein zeta function $\zeta^{(2)}(s,\Delta)$ to the
complex $s$-plane (Section \ref{Sec21}) and basic equations for
determining critical and off-critical zeros (Section \ref{Sec22}).
Section \ref{Sec3} deals with critical zeros of $\zeta^{(2)}(s,\Delta)$.
Based on numerical calculation of open and closed curves of critical
zeros in the $(\Delta,\rho_y)$ plane, the critical edge zeros are introduced
in Section \ref{Sec31}, together with an explicit form of two coupled
integral equations determining their position in the $(\Delta,\rho_y)$ plane.
The singular expansion of $\rho(y)$ around the edge points is discussed in
Section \ref{Sec32}.
An accurate approximation of $\zeta^{(2)}(s,\Delta)$ for small (and large)
values of $\Delta$, indicating an equidistant distribution of critical zeros
along the imaginary axis, is presented in Section \ref{Sec33}.
Section \ref{Sec4} is about off-critical zeros of $\zeta^{(2)}(s,\Delta)$.
The most important result of this work, the generation mechanism of curves of
off-critical zeros starting from critical edge zeros, is explained in
Section \ref{Sec41}.
Pairs of real off-critical zero are found in Section \ref{Sec42}.
    
\renewcommand{\theequation}{2.\arabic{equation}}
\setcounter{equation}{0}

\section{Technicalities} \label{Sec2}

\subsection{Regularization of $\zeta^{(2)}(s,\Delta)$} \label{Sec21}
This part is devoted to a short review of the analytic continuation
of the Epstein zeta function (\ref{sum1}), defined when the real part
$\Re(s)>1$, to the whole complex $s$-plane.
We have the following theorem.

\begin{thm}[\textbf{Analytic continuation of the Epstein zeta function}
\cite{Blanc15,Epstein03}]
Let $\Delta>0$ and $s\in \mathbb{C}$. Let us define
\begin{equation}\label{def:ZsDelta}
Z(s,\Delta) := 
\left( \frac{\Delta}{\pi} \right)^s \Gamma(s) \zeta^{(2)}(s,\Delta) .
\end{equation}

We have
\begin{equation} \label{final}
Z(s,\Delta) = - \frac{1}{2} \left( \frac{1}{1-s} + \frac{1}{s} \right)
+ \frac{1}{2} \int_0^1  \left( t^{s-1} + t^{-s} \right)
\left[ \theta_3\left({\rm e}^{-\pi t\Delta}\right)
\theta_3\left({\rm e}^{-\pi t/\Delta}\right) - \frac{1}{t} \right]{\rm d}t, 
\end{equation} 
where the Jacobi elliptic theta function (with zero argument, see
{\rm \cite{Gradshteyn}}) is given by
\begin{equation} \label{theta3}
\theta_3(q) = \sum_{j\in \mathbb{Z}} q^{j^2} .
\end{equation}
In particular, {\rm (\ref{final})} represents an analytic continuation of
$\zeta^{(2)}(s,\Delta)$ to the whole complex plane, except for
the simple pole at $s=1$.
$Z(s,\Delta)$ fulfills the following duality relation 
\begin{equation} \label{duality:ZsDelta}
Z(s,\Delta) = Z(1-s,\Delta).
\end{equation}
\end{thm}

\begin{remark}
Note that the symmetry 
\begin{equation}
\label{Zsymmetry}
Z(s,\Delta) = Z(s,1/\Delta) 
\end{equation}
is automatically ensured by formula {\rm (\ref{final})}.
\end{remark}

\subsection{Integral equations determining zeros} \label{Sec22}

We start by defining what is a zero of $s\mapsto \zeta^{(2)}(s,\Delta)$
associated to the parameter $\Delta$.

\begin{defi}
We say that $\rho\in \mathbb{C}$ is:
\begin{itemize}
\item a critical zero (of $s\mapsto \zeta^{(2)}(s,\Delta)$) associated to the parameter $\Delta$ if $\rho=\frac{1}{2}+i\rho_y$, $\rho_y\in \mathbb{R}$ and $\zeta^{(2)}(\rho,\Delta)=0$;
\item an off-critical zero (of $s\mapsto \zeta^{(2)}(s,\Delta)$) associated to
the parameter $\Delta$ if $\rho=\rho_x+i\rho_y$, $(\rho_x,\rho_y)\in \mathbb{R}\times \mathbb{R}$, $\rho_x\ne \frac{1}{2}$ and
$\zeta^{(2)}(\rho,\Delta)=0$.
\end{itemize}
\end{defi}
The nontrivial zeros of $\zeta^{(2)}(s,\Delta)$ are related to the
nullity condition of the right-hand side of (\ref{final}). In the case of critical zeros $\rho=\frac{1}{2} + i \rho_y$, the right-hand side of
(\ref{final}) is real and we obtain the following result.

\begin{thm}[\textbf{Equation for critical zeros}]
Let $\rho=\frac{1}{2} + i \rho_y$ be a critical zero associated to the parameter $\Delta$, then $\rho_y$ satisfies the following equation:
\begin{equation} \label{critzeros}
- \frac{2}{1+4\rho_y^2} + \int_0^1 
\cos\left( \rho_y\log t\right) \left[   
\theta_3\left({\rm e}^{-\pi t\Delta}\right)
\theta_3\left({\rm e}^{-\pi t/\Delta}\right) - \frac{1}{t} \right]\frac{{\rm d}t}{\sqrt{t}} = 0 .
\end{equation}  
\end{thm}
\begin{remark}[\textbf{Symmetries of the equation and consequences}]
The symmetry of \eqref{critzeros} with respect to the transformation
$\Delta\to 1/\Delta$ tells us that the set of critical zeros
is the same for the couple of values $\Delta$ and $1/\Delta$.
The symmetry $\rho_y\to -\rho_y$ means that to each critical zero
$\left( \frac{1}{2},\rho_y\right)$ there exists a complex conjugate
critical zero $\left( \frac{1}{2},-\rho_y\right)$.
\end{remark}
In the case of off-critical zeros with $\rho_x\ne \frac{1}{2}$, it is useful
to introduce the deviation of $\rho_x$ from its critical value
\begin{equation}
\delta\rho_x := \rho_x - \frac{1}{2} .
\end{equation}
The right-hand side of (\ref{final}) becomes complex in this case and the
off-critical zeros are determined by a pair of coupled integral equations
(corresponding to the real and imaginary parts).

\begin{thm}[\textbf{Equation for the off-critical zeros}]
Let $\rho=\frac{1}{2}+\delta \rho_x+i\rho_y$, $\delta \rho_x\ne 0$, be an
off-critical zero associated to the parameter $\Delta$, then $\rho$ satisfies
the following equations:
\begin{eqnarray}
- \left[ \frac{1+2\delta\rho_x}{(1+2\delta\rho_x)^2+4\rho_y^2}
+ \frac{1-2\delta\rho_x}{(1-2\delta\rho_x)^2+4\rho_y^2} \right]    
\phantom{aaaaaaaaaaaaaaaaaaa} & & \nonumber \\
+ \int_0^1 
\cos\left( \rho_y\log t\right) \cosh\left( \delta\rho_x\log t\right)
\left[ \theta_3\left({\rm e}^{-\pi t\Delta}\right)
\theta_3\left({\rm e}^{-\pi t/\Delta}\right) - \frac{1}{t} \right]\frac{{\rm d}t}{\sqrt{t}} & = & 0 ,
\label{offzeros1}
\end{eqnarray}
\begin{eqnarray}
2\rho_y \left[ \frac{1}{(1+2\delta\rho_x)^2+4\rho_y^2}
- \frac{1}{(1-2\delta\rho_x)^2+4\rho_y^2} \right]    
\phantom{aaaaaaaaaaaaaaaaa} & & \nonumber \\
+ \int_0^1 
\sin\left( \rho_y\log t\right) \sinh\left( \delta\rho_x\log t\right)
\left[ \theta_3\left({\rm e}^{-\pi t\Delta}\right)
\theta_3\left({\rm e}^{-\pi t/\Delta}\right) - \frac{1}{t} \right]\frac{{\rm d}t}{\sqrt{t}} & = & 0 .
\label{offzeros2}
\end{eqnarray}  
Furthermore, these equations are invariant with respect to the transformations
$\Delta\to 1/\Delta$, $\rho_y\to -\rho_y$ and $\delta\rho_x\to -\delta\rho_x$.
In particular, to each off-critical zero
$\left( \frac{1}{2}+\delta\rho_x,\rho_y\right)$ there exist
the related off-critical zeros
$\left( \frac{1}{2}+\delta\rho_x,-\rho_y\right)$,
$\left( \frac{1}{2}-\delta\rho_x,\rho_y\right)$ and
$\left( \frac{1}{2}-\delta\rho_x,-\rho_y\right)$.
\end{thm}

\renewcommand{\theequation}{3.\arabic{equation}}
\setcounter{equation}{0}

\section{Zeros on the critical line} \label{Sec3}

\subsection{Critical edge zeros} \label{Sec31}
Let us first comment on Figure \ref{allcritintro}. The critical zeros of the Epstein zeta function $\zeta^{(2)}(s,\Delta)$,
calculated numerically by using Eq. (\ref{critzeros}), are represented by
open symbols in Figure \ref{allcritintro}.
The set of zeros lying on the same (closed or open) curve are denoted by
a common open symbol (circle, square, triangle,$\ldots$).
There is an infinite sequence of loop circuits in the left-down corner
of the figure with the common point at the origin $(0,0)$ which are not drawn.  
With regard to the symmetry $\Delta\to 1/\Delta$, each critical zero in
the considered interval $\Delta\in (0,1)$ has a counterpart in the
complementary interval $(1,\infty)$.
To maintain high accuracy of the results (8-20 decimal digits), 
only critical zeros with the component $\rho_y$ smaller than 21
are presented.
The numerical evaluation of one critical zero by using \emph{Mathematica}
takes around 5 seconds of CPU time on a conventional PC.

Each curve of critical zeros in Figure \ref{allcritintro} is defined by
the function $\rho_y(\Delta)$.
Varying the parameter $\Delta$ in the positive or negative direction,
the distance between a couple of nearest zeros can go to zero
and the zeros merge at points referred to as the critical ``edge'' zeros.
They originate at specific values of $\Delta=\Delta^*$ and have imaginary part
$\rho_y^*=\rho_y(\Delta^*)$. 
More precisely, the edge zeros are defined by a divergent tangent as follows.
\begin{defi}[\textbf{Critical edge zeros}]
We call (critical) edge zero any critical zero $\rho=\frac{1}{2}+i\rho_y$ such that there exists $\Delta^*>0$ satisfying
\begin{equation} \label{diverge}
\left. \frac{\rm d}{{\rm d}\Delta} \rho_y(\Delta) \right\vert_{\Delta=\Delta^*}
= \pm \infty .
\end{equation}
Furthermore, assuming that the curve of critical zeros is defined by
the inverse function $\Delta(\rho_y)$, the condition for the edge zeros
{\rm (\ref{diverge})} is equivalent to
\begin{equation} \label{zero}
\left. \frac{\rm d}{{\rm d}\rho_y} \Delta(\rho_y) \right\vert_{\rho_y=\rho_y^*}
= 0
\end{equation}  
and $\Delta^*=\Delta(\rho_y^*)$.
\end{defi}

The edge zeros split into two groups: the left/right edge zeros
are situated on the left/right with respect to the curve $\rho_y(\Delta)$.
More precisely:
\begin{defi}[\textbf{Left/right critical edge zeros}]
Let $\rho^*=\frac{1}{2}+i\rho_y^*$ be a (critical) edge zero
associated to the parameter $\Delta^*$ and the curve of critical
zeros is defined by the inverse function $\Delta(\rho_y)$.
We say that $\rho$ is a right (resp. left) edge zero if
$$
\frac{d^2}{d\rho_y^2} \Delta(\rho_y)\vert_{\rho_y=\rho_y^*}<0
\quad \left(\textnormal{resp. } \frac{d^2}{d\rho_y^2}
\Delta(\rho_y)\vert_{\rho_y=\rho_y^*}>0\right).
$$
\end{defi}

The right and left edge zeros are denoted in Figure \ref{allcritintro} as 
2a, 3a, 4a,$\ldots$ and 1, 2b, 3b, 4b,$\ldots$, respectively;
the close connection between the right and left edge zeros 2a-2b, 3a-3b, etc.
will become clear later (see Section \ref{Sec41}).
Notice that the critical zero numbered 1 with coordinates
$(\Delta^*,\rho_y^*)=(\Delta_c^*\approx 0.141733,0)$ is a left edge zero
because the corresponding curve continues reflection-symmetrically across
the $\Delta$-axis into the lower quadrant. 

We can easily deduce closed-form equations for specifying edge zeros.

\begin{thm}[\textbf{Equation satisfied by a critical edge zero}]
If  $\rho^*=\frac{1}{2}+i\rho_y^*$ is an edge zero associated to
the parameter $\Delta^*$, then
\begin{equation} \label{edge22}
f(\rho_y^*,\Delta^*) = 0 ,
\end{equation}
where the function $f$ is given by
\begin{equation} \label{edge21}
f(\rho_y,\Delta) := \frac{16\rho_y}{\left( 1+4{\rho_y}^2\right)^2} -
\int_0^1  \log t
\sin\left( \rho_y\log t\right) \left[   
\theta_3\left({\rm e}^{-\pi t\Delta}\right)
\theta_3\left({\rm e}^{-\pi t/\Delta}\right) - \frac{1}{t} \right]\frac{{\rm d}t}{\sqrt{t}}.
\end{equation}
\end{thm}
\begin{proof}
We first recall that an edge zero has to satisfy the general
equation (\ref{critzeros}) for critical zeros, i.e.,
\begin{equation} \label{edge1}
- \frac{2}{1+4{\rho_y^*}^2} + \int_0^1 
\cos\left( \rho_y^*\log t\right) \left[   
\theta_3\left({\rm e}^{-\pi t\Delta^*}\right)
\theta_3\left({\rm e}^{-\pi t/\Delta^*}\right) - \frac{1}{t} \right]\frac{{\rm d}t}{\sqrt{t}} = 0 .
\end{equation}
Taking into account the edge-zero condition (\ref{zero}), the derivative of
(\ref{critzeros}) with respect to $\rho_y$ leads to the result.
\end{proof}

It has to be noticed that the coupled equations (\ref{edge22}) and
(\ref{edge1}) have an infinite number of real solutions
for $\Delta^*$ and $\rho_y^*$.
The characteristics of edge critical zeros from Figure \ref{allcritintro},
constrained to the intervals $0<\Delta\leq 1$ and $0\leq \rho_y \leq21$,
are summarized in Table \ref{edgetable}.

\begin{table}[t]
\caption{The coordinates of edge critical zeros appearing in the region of
$0<\Delta\le 1$ and $0\le \rho_y\le 21$ in Figure \ref{allcritintro}.}
\label{edgetable}
\begin{tabular}{cccc}
\hline
edge point &$\Delta^*$&$\rho^*_y$\\
\hline
1&0.141733239663887&0\\
2a&0.798382429865856&9.17479405815734\\
2b&0.950672823506692&13.5092488680816\\
3a&0.578095740200051&11.2961629757333\\
3b&0.688797339793161&12.7134082666419\\
4a&0.708261915413478&14.3461052173020\\
4b&0.810471985748564&17.3035168808027\\
5a&0.577833206956181&15.7904269230734\\
5b&0.625830051933379&16.7721421891791\\
6a&0.397042034784957&14.9386821841068\\
6b&0.499955572107973&16.3629327845743\\
7a&0.634086781531453&17.8588321271621\\
7b&0.690295752437308&19.4462462865857\\
8a&0.356573014664413&16.4816098051657\\
8b&0.375454386384881&16.8495675287149\\
9a&0.337272867689201&18.5674591768417\\
9b&0.374296061779980&19.3629136770424\\
10a&0.560652822542094&19.8540419510498\\
10b&0.577320038404815&20.4237238736290\\
11a&0.896821462590355&19.0008766867965\\
12a&0.778481639573212&20.0611304186419\\
12b&0.789270563104711&20.2667094854061\\
13a&0.578437965650995&20.8178435639014\\
14a&0.309679721075915&20.1102459521285\\
\hline
\end{tabular}
\vspace*{-4pt}
\end{table}

\subsection{Singular expansion around critical edge zeros} \label{Sec32}
The map $\Delta\mapsto \rho_y(\Delta)$ exhibits locally an analytic expansion in $\Delta$,
except for the edge zeros where it shows a singular expansion in $\Delta$.
The singular expansion around edge zeros can be documented by performing
the Taylor series expansion of the inverse function $\Delta(\rho_y)$
on the corresponding curve of critical zeros.
We obtain the following result.

\begin{lemma}[\textbf{First order asymptotics around a critical edge zero}]\label{asympt:edge1storder}
Let $\rho^*=\frac{1}{2}+i\rho_y^*$ be an edge zero associates to the parameter $\Delta^*$. Then, as $\rho_y\to \rho_y^*$ where $\rho=\frac{1}{2}+i\rho_y$ is a critical zero associated to the parameter $\Delta$, we have
$$
\rho_y-\rho_y^* =\mathcal{O}\left( \sqrt{\vert \Delta-\Delta^*\vert}\right).
$$
\end{lemma}
\begin{proof}
The result directly follows from the order two Taylor expansion of the inverse function $\Delta(\rho_y)$ given by
\begin{equation}
\Delta(\rho_y) = \Delta(\rho_y^*) + \frac{1}{2!} \left.
\frac{\rm d^2}{{\rm d}\rho_y^2} \Delta(\rho_y) \right\vert_{\rho_y=\rho_y^*}
\left( \rho_y-\rho_y^* \right)^2 + o\left(\left( \rho_y-\rho_y^* \right)^2 \right),
\end{equation}  
where the condition (\ref{zero}) was taken into account.
\end{proof}

We now perform a general analysis for critical zeros, deriving an asymptotic expansion of Equation (\ref{critzeros}) around a critical zero $\rho=\frac{1}{2}+i\rho_y$ associated to the parameter $\Delta$.

\begin{lemma}[\textbf{Asymptotic expansion around a general critical zero}]
Let $\rho=\frac{1}{2}+i\rho_y$ and $\tilde{\rho}=\frac{1}{2}+i(\rho_y+\delta\rho_y)$ be two critical zeros respectively associated to the parameters $\Delta$ and $\Delta+\delta\Delta$. Then we have, as $\delta\Delta\to 0$ (and therefore $\delta\rho_y\to 0$),
\begin{equation} \label{finalequation}
a \delta\Delta + f \delta\rho_y + c \left( \delta\rho_y \right)^2
- b (\delta\Delta) (\delta\rho_y) - d \left( \delta\rho_y \right)^3
+ {\cal O}\left[\left( \delta\Delta \right)^2 \right]
+ {\cal O}\left[ \delta\Delta\left( \delta\rho_y \right)^2 \right] = 0 ,
\end{equation}
where the function $f=f(\rho_y,\Delta)$ is defined by {\rm (\ref{edge21})}
and the other prefactor functions $a$, $b$, $c$ and $d$,
depending on $(\rho_y,\Delta)$, are given by 
\begin{equation} \label{alpha}
a = \pi \int_0^1\ \sqrt{t} \cos\left( \rho_y \log t \right)
\left[ \theta_3\left( {\rm e}^{-\pi t\Delta} \right) \frac{1}{{\Delta}^2}
\vartheta\left( {\rm e}^{-\pi t/\Delta} \right) -
\theta_3\left( {\rm e}^{-\pi t/\Delta} \right)
\vartheta\left( {\rm e}^{-\pi t\Delta} \right) \right] {\rm d}t , 
\end{equation}  
\begin{equation} \label{beta}
b = \pi \int_0^1  \sqrt{t} \sin\left( \rho_y\log t \right)
(\log t)\left[ \theta_3\left( {\rm e}^{-\pi t\Delta} \right) \frac{1}{{\Delta}^2}
\vartheta\left( {\rm e}^{-\pi t/\Delta} \right) -
\theta_3\left( {\rm e}^{-\pi t/\Delta} \right)
\vartheta\left( {\rm e}^{-\pi t\Delta} \right) \right]{\rm d}t , 
\end{equation}
\begin{equation} \label{gamma}
c = \frac{8\left( 1 - 12 {\rho_y}^2\right)}{
\left(1+4{\rho_y}^2\right)^3} - \frac{1}{2} \int_0^1 
\cos\left( \rho_y\log t\right) \left( \log t\right)^2
\left[ \theta_3\left({\rm e}^{-\pi t\Delta}\right)
\theta_3\left({\rm e}^{-\pi t/\Delta}\right) - \frac{1}{t} \right]
\frac{{\rm d}t}{\sqrt{t}} ,
\end{equation}
\begin{equation} \label{delta}
d = \frac{128\rho_y\left( 1 - 4 {\rho_y}^2\right)}{
\left(1+4{\rho_y}^2\right)^4} - \frac{1}{6} \int_0^1 
\sin\left( \rho_y\log t\right) \left( \log t\right)^3
\left[ \theta_3\left({\rm e}^{-\pi t\Delta}\right)
\theta_3\left({\rm e}^{-\pi t/\Delta}\right) - \frac{1}{t} \right]
\frac{{\rm d}t}{\sqrt{t}}
\end{equation}
with $\theta_3$ is defined by \eqref{theta3} and
\begin{equation} \label{vartheta}
\vartheta(q) := q \frac{{\rm d}}{{\rm d}q} \theta_3(q)
= \sum_{j\in \mathbb{Z}} j^2 q^{j^2}.
\end{equation}
\end{lemma}
\begin{proof}
One substitutes
\begin{equation} \label{infinitesimal}
\Delta \to \Delta + \delta\Delta , \quad \textnormal{and}\qquad
\rho_y \to \rho_y + \delta\rho_y 
\end{equation}
into equation (\ref{critzeros}) and expands in Taylor series in powers of small
deviations $\delta\Delta$ and $\delta\rho_y$.
Furthermore, the $\theta_3$-functions appearing in (\ref{critzeros})
are expanded as follows
\begin{equation}
\theta_3\left( {\rm e}^{-\pi t\Delta} \right) \to
\theta_3\left( {\rm e}^{-\pi t\Delta} \right)
- \pi t (\delta\Delta) \vartheta\left( {\rm e}^{-\pi t\Delta} \right) ,
\quad
\theta_3\left( {\rm e}^{-\pi t/\Delta} \right) \to
\theta_3\left( {\rm e}^{-\pi t/\Delta} \right)
+ \frac{\pi t \delta\Delta}{\Delta^2}
\vartheta\left( {\rm e}^{-\pi t/\Delta} \right) . 
\end{equation}  
It is now straightforward to obtain the desired result.
\end{proof}

It is easy to derive the following asymptotics for $\delta\rho_y$
in the case of a critical zero which does not belong to the set of edge zeros.

\begin{lemma}[\textbf{Asymptotic expansion around a critical non-edge zero}]
Let $\rho=\frac{1}{2}+i(\rho_y+\delta\rho_y)$ be a critical zero associated
to the parameter $\Delta+\delta\Delta$ such that the critical zero
$\frac{1}{2}+i\rho_y$ associated to the parameter $\Delta$ is not an edge zero.
Then we have, as $\delta \Delta\to 0$ (and consequently
$\delta \rho_y\to 0$),
$$
\delta \rho_y = - \frac{a}{f}\delta \Delta +
{\cal O}\left( \delta \rho_y\delta\Delta\right).
$$
\end{lemma}
\begin{proof}
Since $\rho$ is not an edge zero, we have $f(\rho_y,\Delta)\ne 0$ and
therefore the leading term of the expansion of the deviation $\delta\rho_y$
in small $\delta\Delta$ is determined by the equation
$a \delta\Delta + f \delta\rho_y
+ {\cal O}\left(\delta \rho_y\delta\Delta \right)=0$,
which leads to our result.
\end{proof}

We can also derive from (\ref{finalequation}) the following asymptotics around
an edge zero.

\begin{lemma}[\textbf{Asymptotic expansion around a critical edge zero}]
Let $\rho=\frac{1}{2}+i(\rho_y^*+\delta \rho_y)$ associated to the parameter
$\Delta$ be a critical zero  where $\rho^*=\frac{1}{2}+i\rho_y^*$,
associated to the parameter $\Delta^*$ with $\Delta>\Delta^*$, is a left edge
zero.
Then we have, as $\Delta\to \Delta^*$ (and therefore $\rho_y\to \rho_y^*$),
\begin{equation} \label{expansion}
\rho_y(\Delta) = \rho_y(\Delta^*) 
\pm \sqrt{-\frac{a}{c}} \sqrt{\Delta-\Delta^*}
+ \frac{1}{2 c} \left( b - \frac{a d}{c} \right)
\left( \Delta-\Delta^* \right)
+ {\cal O}\left[ \left( \Delta-\Delta^* \right)^{3/2} \right] ,
\end{equation}
where the sign $\pm$ determines the up/down branches of $\rho_y(\Delta)$.
\end{lemma}
\begin{remark}
In particular, the leading order term of our asymptotics is given, as
$\Delta\to \Delta^*$, by
\begin{equation}\label{deltarho1}
\rho_y(\Delta) = \rho_y(\Delta^*) 
\pm \sqrt{-\frac{a}{c}} \sqrt{\Delta-\Delta^*}
+ o\left( \sqrt{\Delta-\Delta^*}\right).
\end{equation}
Here, the ratio $-a/c$ must be positive for the real component $\rho_y$
to exist and it was checked numerically for all left edge zeros
presented in Table \ref{edgetable} that it is so.
Notice that singular expansion of type {\rm (\ref{deltarho1})} with the critical
exponent $\frac{1}{2}$ occurs in a mean-field description of classical
statistical systems at the second-order phase transition
\cite{Baxter82,Samaj13}.
An analogous analysis can be done for right edge zeros.
\end{remark}
\begin{proof}
For a critical edge zero $\rho^*=\frac{1}{2}+i\rho_y^*$, the linear term
of order $\delta\rho_y$ is absent in (\ref{finalequation})
since $f(\rho_y^*,\Delta^*)=0$.
Terms on the left-hand side of (\ref{finalequation}) can be classified
according to their power in the smallness parameter $\delta\Delta$ and
the leading order can be easily derived.
To go to the next order in $\delta\Delta$, one adds to $\delta\rho_y$
in (\ref{deltarho1}) a higher-order term $\alpha \delta\Delta$ with the
as-yet-undetermined constant $\alpha$. 
Expanding all functions in (\ref{finalequation}) up to the order
$\left( \delta\Delta \right)^{3/2}$ one ends up with the relation
$2c \alpha - b + a d/c = 0$ which fixes $\alpha$.
It is now straightforward to get our final expression of the asymptotics.
\end{proof}

\textbf{Application and comparison with the direct computations.} Let us choose one of the left edge critical zeros in Figure \ref{allcritintro},
say the edge point 3b with characteristics $\Delta^*$ and $\rho_y^*$
listed in Table \ref{edgetable}.
After the numerical evaluation of the coefficients $a, b, c$ and $d$,
the expansion (\ref{expansion}) takes the form
\begin{equation} \label{expansion3b}
\delta\rho_y(\Delta) = 
\pm 4.87411 \sqrt{\Delta-\Delta^*} + 
22.493 \left( \Delta-\Delta^* \right)
+ {\cal O}\left[ \left( \Delta-\Delta^* \right)^{3/2} \right] . 
\end{equation}
The fitting of numerical data for the up branch (the $+$ sign) implies
the prefactor $4.87412$ to $\sqrt{\Delta-\Delta^*}$ and the one
$22.686$ to $\left( \Delta-\Delta^* \right)$, while for the down branch
(the $-$ sign) the corresponding prefactors $-4.87415$ and $22.230$
deviate a bit more from the ``exact'' ones in (\ref{expansion3b}),
but they are still suitable. 

\subsection{Approximation of Epstein zeta function for small values of
$\Delta$} \label{Sec33}
In this part, for small values of the anisotropy parameter $\Delta$,
the Epstein function is approximated well by a function which reveals
an equidistant distribution of critical zeros along the imaginary axis
in the limit $\Delta\to 0$. We start by showing the following asymptotic result as $\Delta\to 0$.

\begin{thm}[\textbf{Asymptotics for small values of $\Delta$}]
As $\Delta\to 0$, we have
\begin{equation} \label{deltasmall1}
\zeta^{(2)}(s,\Delta) = \frac{1}{\Delta^{2s}} \frac{1}{\pi^{\frac{1}{2}-2s}}
\frac{\Gamma\left( \frac{1}{2}-s \right)}{\Gamma(s)} \zeta(1-2s) 
+ \frac{1}{\Delta} 
\frac{\sqrt{\pi} \Gamma\left( s-\frac{1}{2}\right)}{\Gamma(s)} \zeta(2s-1)
+ \mathcal{O}(\Delta). 
\end{equation}
\end{thm}
\begin{proof}
The sum in (\ref{sum1}) can be straightforwardly converted to
\begin{equation}
\frac{1}{2} \xsum_{j,k=-\infty}^{\infty} \frac{1}{(j^2+\Delta^2 k^2)^s} =
\zeta(2s) + \frac{1}{\Delta^{2s}} \zeta(2s) +
2 \sum_{j,k=1}^{\infty} \frac{1}{(j^2+\Delta^2 k^2)^s} .
\end{equation}  
Writing
\begin{equation}
\sum_{k=1}^{\infty} \frac{1}{(j^2+\Delta^2 k^2)^s} =
\sum_{k=0}^{\infty} \frac{1}{(j^2+\Delta^2 k^2)^s} - \frac{1}{j^{2s}} 
\end{equation}
the Epstein zeta function is expressed as
\begin{equation}
\zeta^{(2)}(s,\Delta) = \frac{1}{\Delta^{2s}} \zeta(2s) - \zeta(2s)
+ \sum_{j=1}^{\infty} \sum_{k=0}^{\infty} f_j(k) ,
\end{equation}    
where
\begin{equation} \label{fj}
f_j(x)  =  \frac{2}{(j^2+\Delta^2 x^2)^s} .  
\end{equation}  
In the limit $\Delta\to 0$, the difference between successive terms
in the sum $\sum_{k=0}^{\infty} f_j(k)$ is also negligibly small.
The sum can be thus treated as an integral according to
the Euler-Maclaurin formula \cite{Apostol99}
\begin{equation}
\sum_{k=0}^{\infty} f_j(k) = \int_0^{\infty} f_j(x)dx +
\frac{f_j(0)+\displaystyle \lim_{M\to \infty}f_j(M)}{2} + \sum_{l=1}^{\left[\frac{p}{2}\right]}
\frac{B_{2l}}{(2l)!} \left[ \lim_{M \to \infty} f_j^{(2l-1)}(M) - f_j^{(2l-1)}(0) \right] + R_j(p), 
\end{equation}
where $p\geq 2$ is an integer , $\left[ \cdots \right]$ denotes the integer
part, $\{ B_{2l}\}$ are Bernoulli numbers and the absolute value of
the error term $R_j(p)$ is bounded by 
\begin{equation} \label{bound}
\left\vert R_j(p) \right\vert \le \frac{2\zeta(p)}{(2\pi)^p}
\int_0^{\infty} \left\vert f_j^{(p)}(x) \right\vert dx. 
\end{equation}
$f_j(x)$ is an even function of $x$ and therefore its odd derivatives
with respect to $x$ vanish at $x=0$.
Since the derivatives of $f_j(x)$ vanish also at $x\to\infty$, it holds that 
\begin{equation}
\sum_{k=0}^{\infty} f_j(k) = 2 \int_0^{\infty} 
\frac{{\rm d}x}{(j^2+\Delta^2 x^2)^s} + \frac{1}{j^{2s}} + R_j(p) .
\end{equation}
Evaluating the integral
\begin{equation}
\int_0^{\infty} \frac{{\rm d}x}{(j^2+\Delta^2 x^2)^s}
= \frac{1}{\Delta} \frac{1}{j^{2s-1}} \frac{\sqrt{\pi}}{2}
\frac{\Gamma\left( s-\frac{1}{2}\right)}{\Gamma(s)} ,
\end{equation}
we end up with
\begin{equation}
\zeta^{(2)}(s,\Delta) = \frac{1}{\Delta^{2s}} \zeta(2s) + \frac{1}{\Delta}
\frac{\sqrt{\pi} \Gamma\left( s-\frac{1}{2}\right)}{\Gamma(s)} \zeta(2s-1)
+ \sum_{j=1}^{\infty} R_j(p) . 
\end{equation}
By applying the dual relation (see e.g. \cite{Blanc15})
\begin{equation}
\frac{1}{\pi^s} \Gamma(s) \zeta(2s) = \frac{1}{\pi^{\frac{1}{2}-s}}
\Gamma\left( \frac{1}{2}-s \right) \zeta(1-2s) 
\end{equation}
we end up with the following asymptotic expansion as $\Delta\to 0$:
\begin{equation} \label{deltasmall2}
\zeta^{(2)}(s,\Delta) = \frac{1}{\Delta^{2s}} \frac{1}{\pi^{\frac{1}{2}-2s}}
\frac{\Gamma\left( \frac{1}{2}-s \right)}{\Gamma(s)} \zeta(1-2s) 
+ \frac{1}{\Delta} 
\frac{\sqrt{\pi} \Gamma\left( s-\frac{1}{2}\right)}{\Gamma(s)} \zeta(2s-1)
+ \sum_{j=1}^{\infty} R_j(p) . 
\end{equation}
To estimate the error term, writing $s = s_x + i s_y$, we insert $f_j(x)$
(\ref{fj}) into the bound (\ref{bound}) for
$p=2$ and obtain by a straightforward computation that
\begin{equation}
\left|\sum_{j=1}^\infty R_j(2)\right|\leq \sum_{j=1}^\infty
\frac{2\zeta(2)}{(2\pi)^2}\int_0^\infty \left| f_j^{(2)}(x)\right| dx
= \frac{8\Delta\zeta(2)\zeta(2s_x+1)}{(2\pi)^2} \sqrt{s_x^2+s_y^2}  
I(s_x,s_y) ,
\end{equation}
where the integral
\begin{equation}
I(s_x,s_y) = \int_0^1 \left\{ \sqrt{\left[ 1-(1+2 s_x)t^2 \right]^2+4 s_y^2t^4}
+ t^{2 s_x} \sqrt{\left[ t^2-(1+2s_x)\right]^2+4 s_y^2} \right\}
\frac{{\rm d}t}{(1+t^2)^{2+s_x}}
\end{equation}
converges if $s_x>-\frac{1}{2}$.
Note that $\zeta(2 s_x+1)$ diverges for $s_x=0$.
This completes the proof.
\end{proof}

\begin{figure}[t]
\centering
\includegraphics[clip,width=0.75\textwidth]{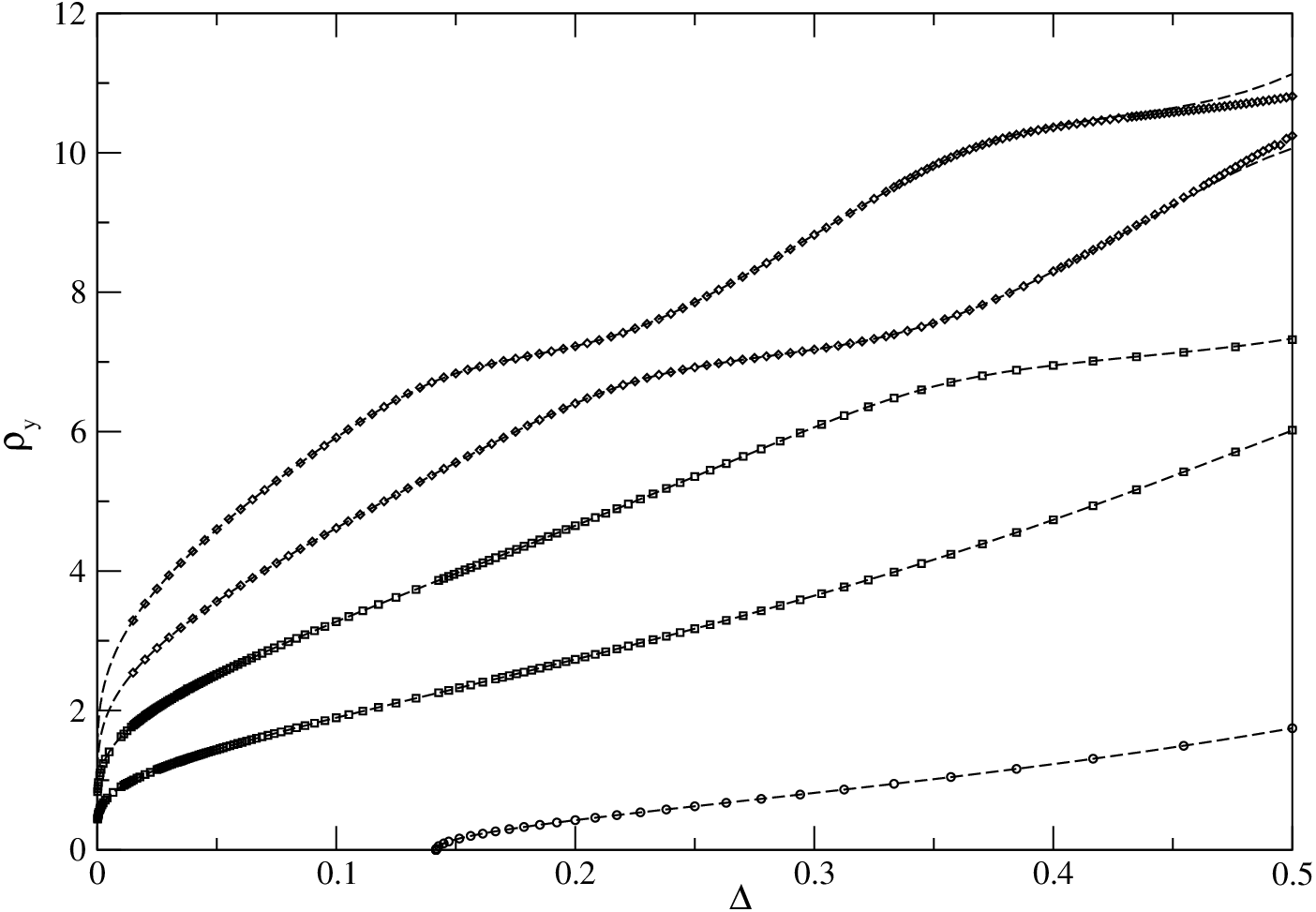}
\caption{The comparison of the results for the first critical zeros as
the functions of $\Delta$ calculated from the approximative equation
(\ref{critappr}) (dashed curves) with the ones obtained by using the exact
equation (\ref{critzeros}) (the symbol notation is taken from
Figure \ref{allcritintro}).}
\label{comparison}
\end{figure}
\begin{remark}
The expansion (\ref{deltasmall2}) then has the meaning of a systematic
Laurent series expansion in $\Delta$ provided that $\Re(s)>0$.
\end{remark}
\textbf{Approximate equation for zeros in the small $\Delta$ regime and comparison with direct computations.}
Let us neglect the error term of order $\Delta$ and consider the approximate
equation for zeros $\{ \rho\}$ of the Epstein zeta function in the region
of small $\Delta$:
\begin{equation} \label{zz}
\left( \frac{\Delta}{\pi} \right)^{2\rho-1} = -
\frac{\Gamma\left( \frac{1}{2}-\rho \right) \zeta(1-2\rho)}{
\Gamma\left( \rho-\frac{1}{2}\right) \zeta(2\rho-1)} .      
\end{equation}  
For the critical zeros $\rho=\frac{1}{2} + i\rho_y$, this equation
takes the form
\begin{equation} \label{critappr}
\left( \frac{\Delta}{\pi} \right)^{2i\rho_y} = -
\frac{\Gamma\left( -i\rho_y \right) \zeta(-2i\rho_y)}{
\Gamma\left( i\rho_y\right) \zeta(2i\rho_y)} .      
\end{equation}
In Figure \ref{comparison}, the results for the first critical zeros as
the functions of the anisotropy parameter $\Delta$ calculated from this
approximative equation (dashed curves) are compared with the ones obtained
by using the exact equation (\ref{critzeros})
(the symbol notation is taken from Figure \ref{allcritintro}).
It is seen that the critical zeros obtained from the approximate equation
(\ref{critappr}) agree with the exact ones unexpectedly far away, up to
$\Delta\approx 0.5$.
Roughly speaking, the approximate formula (\ref{critappr}) works well
until approaching an edge zero at which two critical
zeros merge, the phenomenon which is out of reach of this formula.
This can be seen in the upper right corner of Fig. \ref{comparison}
where the two ``exact'' curves tend to the right edge point 3a. 
High accuracy of the approximative values of the first eight critical zeros
for the Epstein zeta function is documented for
$\Delta=1/\sqrt{7}\approx 0.378$ in Table \ref{approxtable}.
As is intuitively expected, the accuracy of the approximative results
deteriorates as the value of $\rho_y$ increases.

\begin{table}[t]
\caption{The comparison of the exact and approximate values of the first
8 zeros for $\Delta=1/\sqrt{7}$.}
\label{approxtable}
\begin{tabular}{cccc}
\hline
zero \# & $\rho_y\ \ {\rm exact}$&$\rho_y\ \ {\rm approx}$\\
\hline
1&1.133090035457&1.133090358285\\
2&4.475738283729&4.475726461185\\
3&6.845491712491&6.845712742060\\
4&7.931630248198&7.930996972746\\
5&10.19781031911&10.20336832640\\
6&11.16018454312&11.14537554655\\
7&12.48960334303&12.51829228147\\
8&14.13472514173&14.05004856679\\
\hline
\end{tabular}
\vspace*{-4pt}
\end{table}

\medskip

The following result illustrates what we observe on Figure \ref{comparison} and shows that, in the limit $\Delta\to 0$, the set of critical zeros is equidistant along the imaginary axis, with the spacing $\frac{\pi }{\vert \log(\Delta/\pi)\vert}$ between the nearest neighbors.

\begin{lemma}[\textbf{Asymptotic behavior of $\rho_y$ as $\Delta\to 0$}]
Let $\rho=\frac{1}{2}+i\rho_y$ be a critical zero associated to the parameter
$\Delta$.
Then we have 
$$
\lim_{\Delta\to 0} \rho_y(\Delta)=0.
$$
Furthermore, in the limit $\Delta\to 0$, the asymptotic critical zeros are
given by $\{\rho(n):=\frac{1}{2}+\rho_y(n)\}$ where
\begin{equation} \label{sequence}
\rho_y(n) := \frac{\pi n}{\vert \log(\Delta/\pi)\vert} , \qquad
n\in \mathbb{Z}^*.  
\end{equation} 
\end{lemma}

\begin{proof}
Let us consider the two first orders of the asymptotics in \eqref{deltasmall2}
when $s=\frac{1}{2}+i\rho_y$:
\begin{equation}
\zeta^{(2)}(s,\Delta) = \frac{\sqrt{\pi}\Gamma(-i\rho_y)\zeta(-2i\rho_y)}{
\Delta\Gamma\left( \frac{1}{2}+i\rho_y\right)} \left[
\left( \frac{\Delta}{\pi} \right)^{-2i\rho_y} +
\frac{\Gamma(i\rho_y)\zeta(2i\rho_y)}{\Gamma(-i\rho_y)\zeta(-2i\rho_y)}
\right] .
\end{equation}  
Let us assume that in the limit $\Delta\to 0$ also the component $\rho_y$
of critical zeros goes to 0, as is seen in Figure \ref{comparison}.
Since
\begin{equation}
\lim_{\rho_y\to 0} \frac{\Gamma\left( i \rho_y \right) \zeta(2 i \rho_y)}{
\Gamma\left( - i \rho_y\right) \zeta(-2 i \rho_y)} = - 1 ,      
\end{equation}
the critical zeros are given by $(\Delta/\pi)^{2 i\rho_y} = 1$,
in agreement with (\ref{sequence}).
As
\begin{equation}
\frac{\sqrt{\pi}\Gamma(-i\rho_y)\zeta(-2i\rho_y)}{
\Delta\Gamma\left( \frac{1}{2}+i\rho_y\right)} \mathop{\sim}_{\rho_y\to 0}
- \frac{i}{2\rho_y} ,
\end{equation}
the zeta function $\zeta^{(2)}\left( \frac{1}{2}+\rho_y,\Delta\right)$
diverges as $\vert \log\Delta\vert/\Delta$ on the curves of critical
zeros as $\Delta\to 0$, but this has no impact on the location of its zeros.
The proof is complete.
\end{proof}

The distance between the nearest neighbors in
the asymptotic sequence of critical zeros (\ref{sequence}) is
predicted to be $\pi/\vert\log(\Delta/\pi)\vert$; the dependence on
the inverse of $\log(\Delta)$ indicates that one has to take
extremely small values of $\Delta$ to obtain reliable results.
We have performed numerical evaluation of the first
few critical zeros with $\rho_y>0$ for $\Delta=0.0001$ by using
the exact formula (\ref{critzeros}).
The distance between the first and second zeros is 0.375,
between the second and third zeros is 0.357, 
between the third and fourth zeros is 0.347,
between the fourth and fifth zeros is 0.340,
between the fifth and sixth zeros is 0.337, etc.,
which means that the spectrum of zeros is almost equidistant
as was anticipated.
Our asymptotic result (\ref{sequence}) suggests that the distance should be
0.303 which is a reasonable estimate for the considered (not small enough)
value of $\Delta=0.0001$.

\begin{remark}[\textbf{The large $\Delta$ case}]
With regard to the symmetry $\Delta\to 1/\Delta$ of basic equations
{\rm (\ref{critzeros})} for critical and {\rm (\ref{offzeros1})},
{\rm (\ref{offzeros2})} for off-critical zeros, one can accomplish
an analogous analysis in the opposite limit $\Delta\to\infty$, with the result
\begin{equation} \label{zzz}
\left( \pi \Delta \right)^{2\rho-1} = -
\frac{\Gamma\left( \rho-\frac{1}{2} \right) \zeta(2\rho-1)}{
\Gamma\left( \frac{1}{2}-\rho \right) \zeta(1-2\rho)} .      
\end{equation}
Similarly as in the limit $\Delta\to 0$, the spectrum of critical
zeros is equidistant along the imaginary axis in the limit $\Delta\to\infty$,
with the spacing $\frac{\pi}{\log(\pi \Delta)}$ between the nearest neighbors.
\end{remark}

\renewcommand{\theequation}{4.\arabic{equation}}
\setcounter{equation}{0}

\section{Zeros off the critical line} \label{Sec4}

\subsection{Generation of off-critical zeros from critical edge zeros} \label{Sec41}
We observe the following in Figure \ref{allcritintro}: given an edge zero
associated to the parameter $\Delta^*$, there exists $\delta_0>0$ such that
for all $0<\delta<\delta_0$, for a left (resp. right) edge zero, there is no
critical zero associated to the parameter $\Delta^*-\delta$
(resp. $\Delta^*+\delta$).
In the specific case of left edge zeros, this is caused by
the fact that Eq. (\ref{finalequation}) with
the numerically verified inequality $-a/c>0$ has no real
solution for $\delta\rho_y$ if $\delta\Delta = \Delta-\Delta^*<0$, see also
Eq. (\ref{deltarho1}).

Therefore, the only way to have a zero
of $\zeta^{(2)}(s,\Delta)$
corresponding to these values of
the parameter $\Delta$ close to an edge zero is to allow the $\rho_x$-component
to deviate from its critical value $\frac{1}{2}$.
We therefore obtain the following result.

\begin{lemma}[\textbf{Asymptotic expansion of an off-critical zero around a critical zero}]
Let $\rho=\frac{1}{2}+\delta \rho_x+i\left(\rho_y+\delta\rho_y \right)$ be
an off-critical zero associated to the parameter $\Delta+\delta \Delta$ such
that $\tilde{\rho}=\frac{1}{2}+i\rho_y$ is a critical zero associated to
the parameter $\Delta$.
Then we have, as $\delta \Delta\to 0$ (and therefore $\delta \rho_y\to 0$ and
$\delta\rho_x\to 0$), the following two equations:
\begin{equation} \label{first}
a \delta\Delta + f \delta\rho_y - c \left( \delta\rho_x\right)^2
+ c \left( \delta\rho_y\right)^2 - b \delta\Delta \delta\rho_y
+ 3 d \left( \delta\rho_x\right)^2 \delta\rho_y
- d \left( \delta\rho_y\right)^3 + o(\delta\rho_y^3) = 0 ,
\end{equation}  
\begin{equation} \label{second}
\delta\rho_x \left[ -f + b \delta\Delta - 2 c \delta\rho_y
- d \left( \delta\rho_x \right)^2 + 3 d  \left( \delta\rho_y \right)^2
+ o(\delta\rho_y^2)+o(\delta\rho_x^2) \right] = 0 ,    
\end{equation}
where the function $f$ is defined by (\ref{edge21}) and the functions $a$, $b$,
$c$ and $d$ by equations {\rm (\ref{alpha})--(\ref{delta})}.

\noindent Furthermore if $\tilde{\rho}$ is not an edge zero, then there is no such off-critical zero $\rho$ in its neighborhood.
\end{lemma}
\begin{proof}
We simply use (\ref{infinitesimal}) in (\ref{offzeros1})
and (\ref{offzeros2}) and Taylor expanding in powers of small variables $\delta\Delta$, $\delta\rho_x$ and $\delta\rho_y$. Furthermore, if the critical zero $\tilde{\rho}$ is not an edge zero, it holds that
$f(\rho_y,\Delta)\ne 0$.
The second equation (\ref{second}) with sufficiently small $\delta\Delta$
and $\delta\rho_y$ has the only solution $\delta\rho_x=0$.
In other words, there are no off-critical zeros in the neighborhood
of the critical zero which is not an edge zero.
\end{proof}

\begin{lemma}[\textbf{Asymptotic expansion of an off critical zero around a left edge zero}]
Let $\rho=\frac{1}{2}+\delta \rho_x+i\left(\rho_y^*+\delta\rho_y \right)$ be an off-critical zero associated to the parameter $\Delta$ such that $\rho^*=\frac{1}{2}+i\rho_y^*$ is a left edge zero associated to the parameter $\Delta^*$. Then, we obtain, for $\Delta < \Delta^*$, $\Delta\to \Delta^*$ (and then $\delta \rho_y\to 0$ and $\delta \rho_x\to 0$),
\begin{equation} \label{rhox}
\delta\rho_x = \pm \sqrt{-\frac{a}{c}} \sqrt{\Delta^*-\Delta}+o\left(\sqrt{\Delta^*-\Delta} \right),
\end{equation}
and
\begin{equation} \label{rhoy}
\delta\rho_y = - \frac{1}{2c} \left( b - \frac{ad}{c} \right)
\left( \Delta^*-\Delta \right)+ o\left( \Delta^*-\Delta \right) , 
\end{equation} 
\end{lemma}
\begin{remark}
The relations (\ref{rhox}) and (\ref{rhoy}) are the asymptotic formulas
for a curve of off-critical zeros starting from the considered left
edge point which are valid for $\Delta$ close to $\Delta^*$ and $\rho_y$ close
to $\rho_y^*$.

The $\pm$ sign for $\delta\rho_x$ in (\ref{rhox}) means that
at each point along a curve of off-critical zeros $\rho_y(\Delta)$  
there exist a conjugate pair of solutions $\rho_x=\frac{1}{2}+\delta\rho_x$
and $1-\rho_x=\frac{1}{2}-\delta\rho_x$. 
A similar analysis can be made for right edge zeros.
\end{remark}
\begin{proof}
Since $\rho^*$ is an edge zero, we have $f(\rho_y^*,\Delta^*)=0$ and it follows
that (\ref{second}) is satisfied also for $\delta\rho_x\ne 0$.
Indeed, note that in the first equation (\ref{first})
the term $c \left( \delta\rho_y\right)^2$, which was dominant in the previous
analysis of critical zeros in Section \ref{Sec32}, has a counterpart with
the opposite sign $-c \left( \delta\rho_x\right)^2$.
This latter term becomes dominant when the difference
$\delta\Delta=\Delta-\Delta^*$ changes its positive sign to the negative one,
implying that
$$
a(\Delta-\Delta^*) - c(\delta\rho_x)^2 + c(\delta\rho_y)^2
+o(\delta \rho_y^2)+o(\delta \rho_x^2)=0,
$$
Since $\delta\rho_x\ne 0$, the second equation (\ref{second}) implies
that 
$$
b(\Delta-\Delta^*)-2c\delta\rho_y
-d\left( \delta\rho_x \right)^2 + o\left( (\Delta^*-\Delta)^2\right)
+o(\delta \rho_x^2)=0.
$$
These two equations exhibit the solutions with expansions
of type (\ref{rhox}) and (\ref{rhoy}).
\end{proof}

\textbf{Numerical method to generate off-critical zeros curves.}
The fact that each curve of off-critical zeros starts/ends at edge points
simplifies very much the numerical evaluation of off-critical zeros by using
\emph{Mathematica}.
As a function to deal with we take the sum of the squared left-hand side
of Eqs. (\ref{offzeros1}) and (\ref{offzeros2}).
Applying the command FindMinimum to this function, the zero is taken as
sure if the function value is less than $10^{-23}$.
To avoid escape from a local minimum, one starts from a (say right) edge
point and increases $\Delta$ by a tiny amount $0.0001$, after few steps
the shift can be augmented to $0.001-0.01$.
The search for a minimum takes around 60 sec of CPU time on a conventional PC.
For integer values of $\Delta^2=5,6,8,\ldots$, when
off-critical zeros can be calculated with a high precision from
exact sums of products of Dirichlet $L$-functions \cite{Borwein13},
our numerical results agree with these analytic predictions by at least
20 decimal digits.

\medskip

\textbf{Numerical observations.} As is seen in Figure \ref{critoffintro},
each curve of off-critical zeros joins
a pair of critical right (notation ``a'') and left (``b'') edge zeros.
As a rule, the $\rho_y$-coordinate of the right edge point is
smaller than that of the corresponding left edge point.
In the large majority of cases the curves of off-critical zeros go up
monotonously when increasing $\Delta$; the only exception from the curves
presented in Figure \ref{critoffintro} is the curve starting at the right
edge point 11a which first goes down in a short interval of $\Delta$-values
and then goes up to the left edge point 11b (not in the figure).
The intersection of a solid curve of off-critical zeros with a curve of
critical zeros (symbols) is not contradictory: the component $\rho_x$
varies along the solid curves (not indicated in the figure) while it is
constant $\frac{1}{2}$ along the critical curves.
The vertical dashed lines, pictured at the values $1/\sqrt{2}$,
$1/\sqrt{3}$, $1/\sqrt{4}$ and $1/\sqrt{7}$ of $\Delta$, correspond, together
with $\Delta=1$, to the special cases when $\zeta^{(2)}(s,\Delta)$
factorizes itself into product of a zeta function, a Dirichlet $L$ function
and a prefactor function whose zeros (lying on the critical line only)
can be determined analytically \cite{Borwein13,McPhedran16}.
According to the generalized Riemann hypothesis \cite{Davenport36,McPhedran16},
$\zeta^{(2)}(s,\Delta)$ exhibits only critical zeros for these values
of $\Delta$.
This fact is clearly seen in Figure \ref{critoffintro} where no solid curve
of off-critical zeros intersects dashed and $\Delta=1$ lines, although
some of the edge zeros are localized very close to dashed lines.

\begin{figure}[t]
\centering
\includegraphics[clip,width=0.75\textwidth]{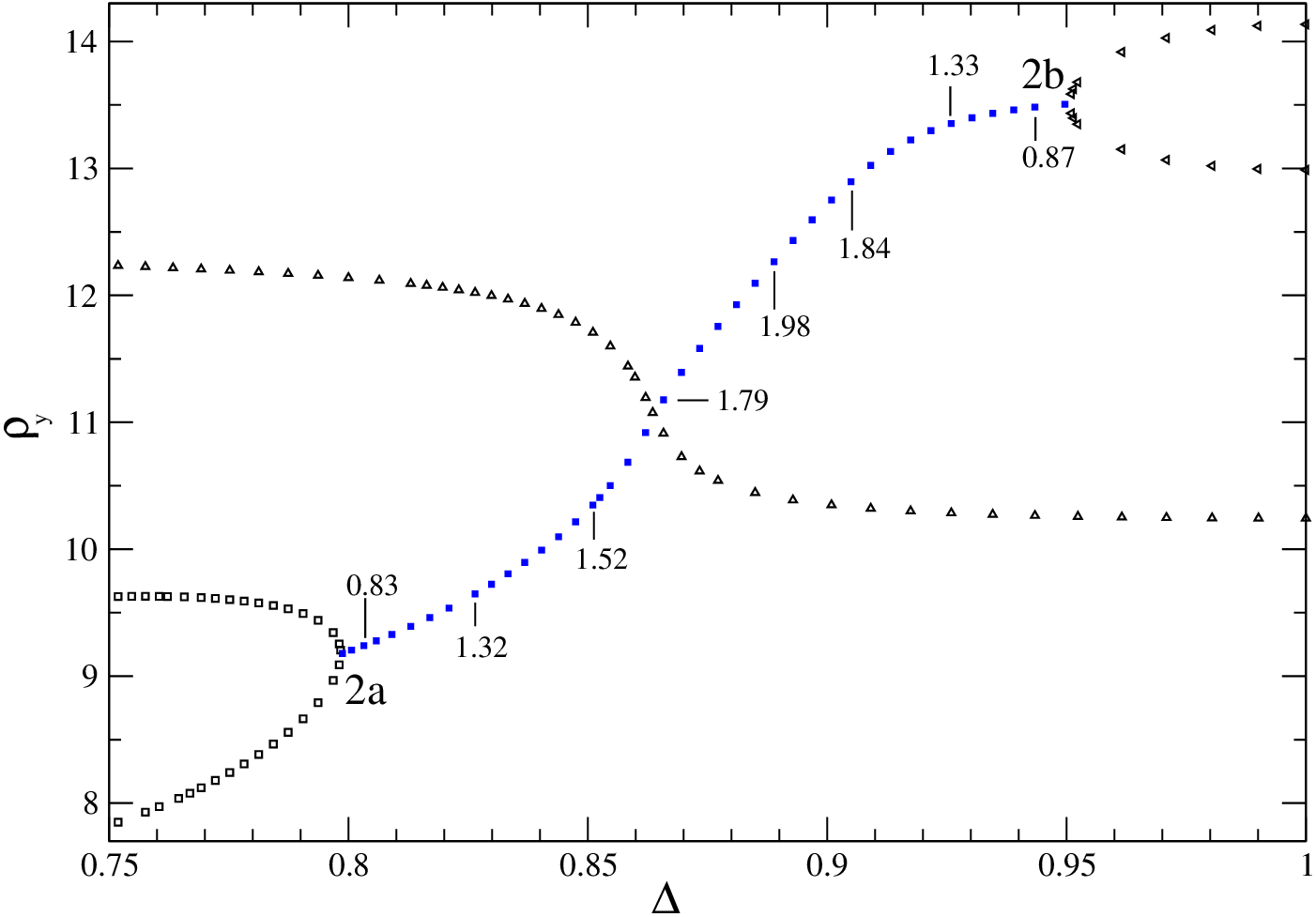}
\caption{Numerical results for the curve off-critical zeros (full squares)
between the right edge point 2a (lying on the curve of critical zeros
represented by open squares) and the left edge point 2b (lying on the curve
of critical zeros represented by open triangles).
The evolution of the component $\rho_x$ along the curve is indicated
by short lines.
At each point of the curve there exists another solution with
the component $1-\rho_x$.}
\label{2a2b}
\end{figure}

Figure \ref{2a2b} documents numerical results for the curve of
off-critical zeros (full squares) interconnecting the pair of right and left
edge zeros, denoted as 2a and 2b in Figure \ref{allcritintro}.
The evolution of the component $\rho_x$ along the curve is indicated
by numbers with short lines attached; by definition of critical zeros,
$\rho_x=\frac{1}{2}$ at the edge points 2a and 2b.
Note that at each point of the curve there exists another solution with
the component $1-\rho_x$.
The intersection of the curve of off-critical zeros (full squares)
with $\rho_x\ne \frac{1}{2}$ and the curve of critical zeros (open triangles)
with $\rho_x=\frac{1}{2}$ is artificial because the corresponding components
$\rho_x$ do not coincide.

\begin{figure}[t]
\centering
\includegraphics[clip,width=0.75\textwidth]{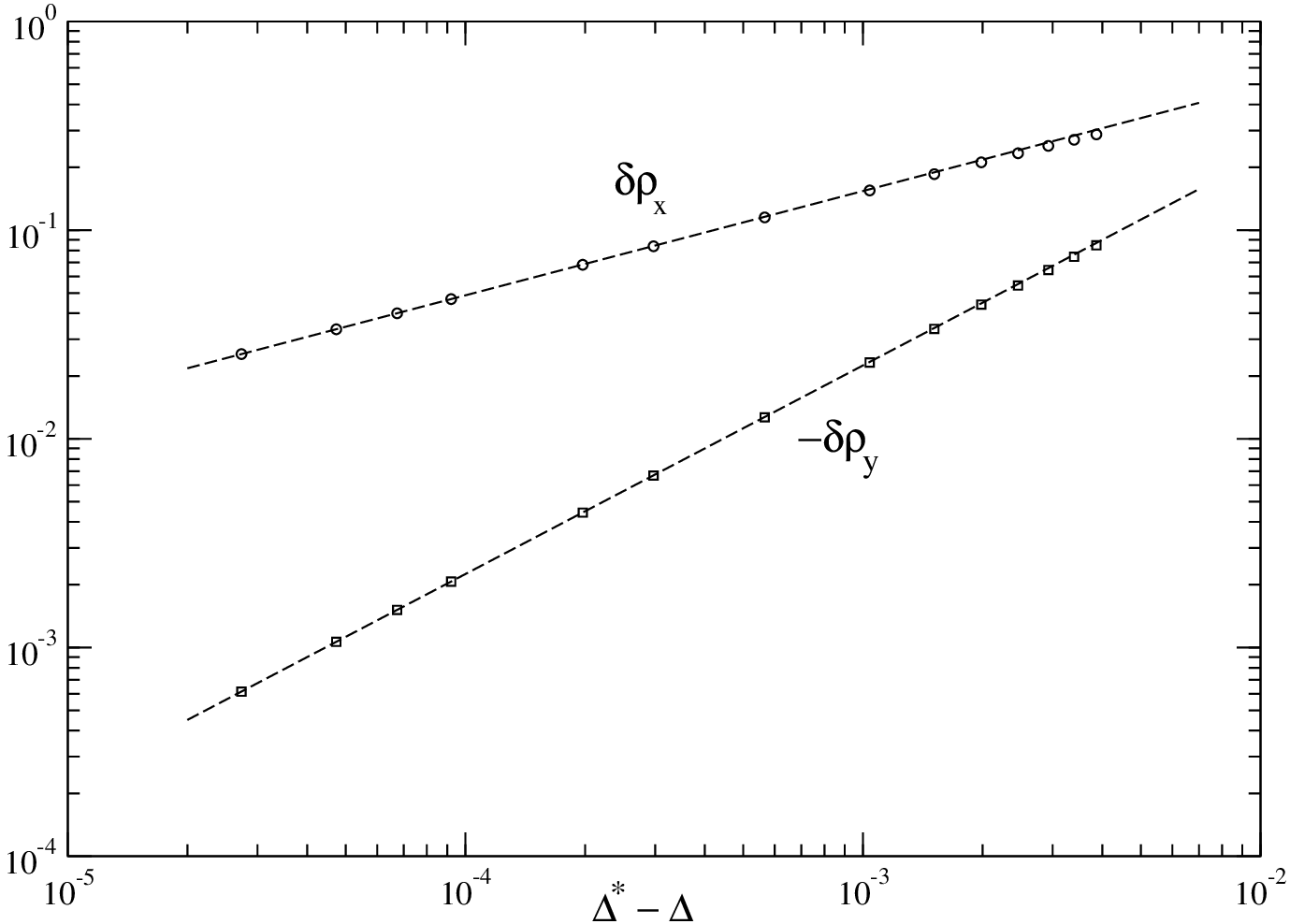}
\caption{The log-log plot of numerical and analytic data for the curve of
off-critical zeros going from the left edge denoted as 3b in Figure
\ref{allcritintro}.  
Numerical dependences of $\delta\rho_x$ (open circles) and $-\delta\rho_y$
(open squares) on the small deviations from the edge point $\Delta^*-\Delta$
are compared with the analytic predictions (\ref{depend}) represented
by dashed lines.} 
\label{num3b}
\end{figure}

\medskip

\textbf{Comparison of our asymptotic formulas and analytic data obtained directly.} A check of the asymptotic formulas (\ref{rhox}) and (\ref{rhoy}) for the curve
of off-critical zeros close to a left edge zero was made for the edge point
denoted as 3b in Figure \ref{allcritintro}, similarly as in the previous case of
the expansion formula (\ref{expansion3b}) for the curve of critical zeros.
The expected dependences
\begin{equation} \label{depend}
\delta\rho_x(\Delta) = \pm 4.87411 \sqrt{\Delta^*-\Delta} , \qquad
-\delta\rho_y(\Delta) = 22.493 \left( \Delta^*-\Delta \right) 
\end{equation}  
are reproduced very well by fitting numerical data,
namely the prefactor obtained for $\delta\rho_x$ equals to $4.87412$ and
the prefactor for $-\delta\rho_y$ is estimated to $22.498$.
The agreement of the asymptotic relations (\ref{depend}) (dashed lines)
with the numerical data (open symbols) is pictured in Figure \ref{num3b}. 

\subsection{Real off-critical zeros} \label{Sec42}
In this section, we are interested in real off-critical zeros. It has already been mentioned that the critical zero numbered by 1 in Figure \ref{allcritintro}, lying on the
$\Delta$-axis, is a left edge zero because
the curve of critical zeros passes across the $\Delta$-axis into
the lower quadrant in a reflection-symmetrical way.

\medskip

\begin{figure}[t]
\centering
\setbox1=\hbox{\includegraphics[clip,height=8.5cm]{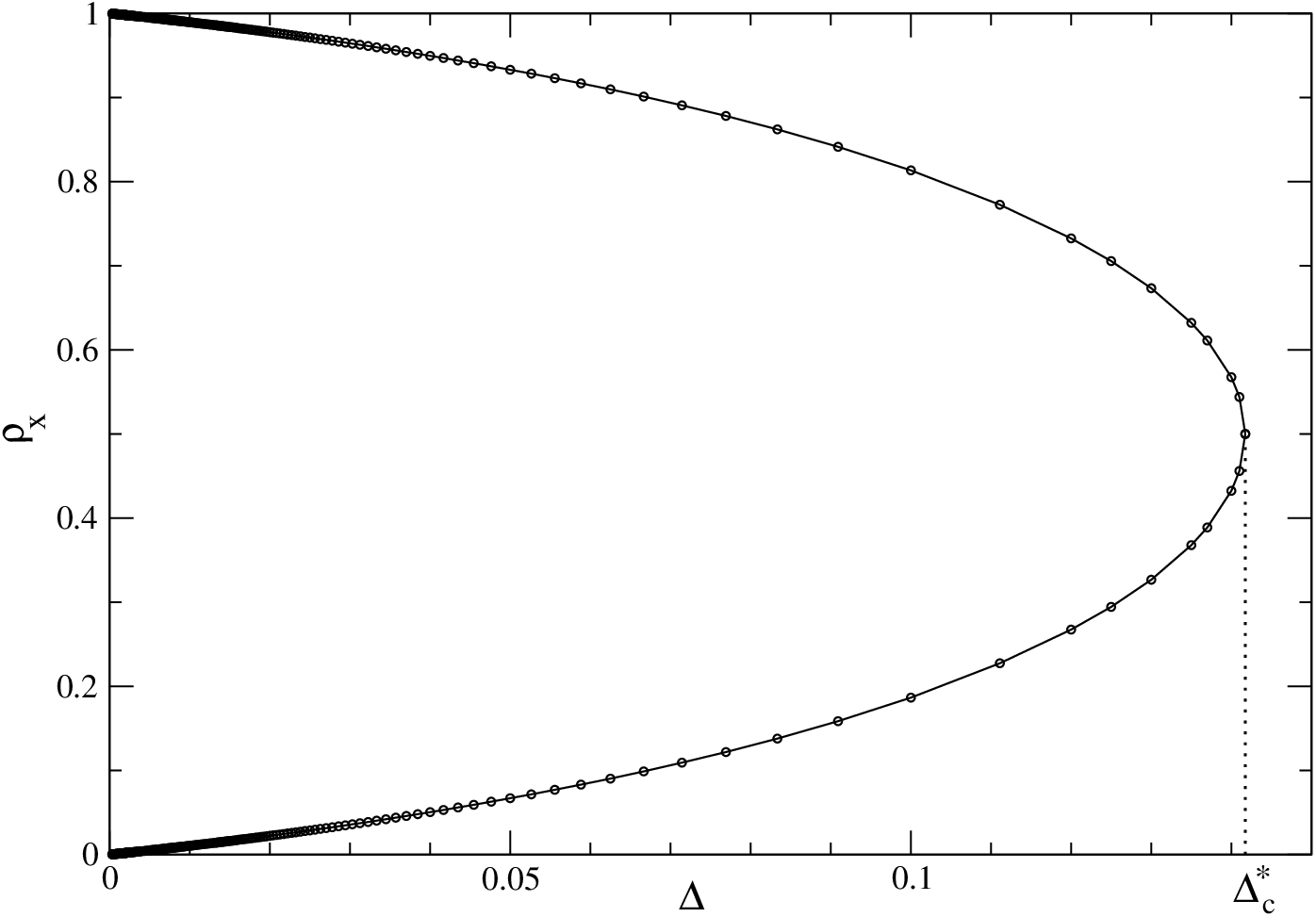}}
\includegraphics[clip,height=8.5cm]{offcritreal.eps}\lapbox[5cm]{-10cm}{
\raisebox{2.3cm}{\includegraphics[clip,height=4.1cm]{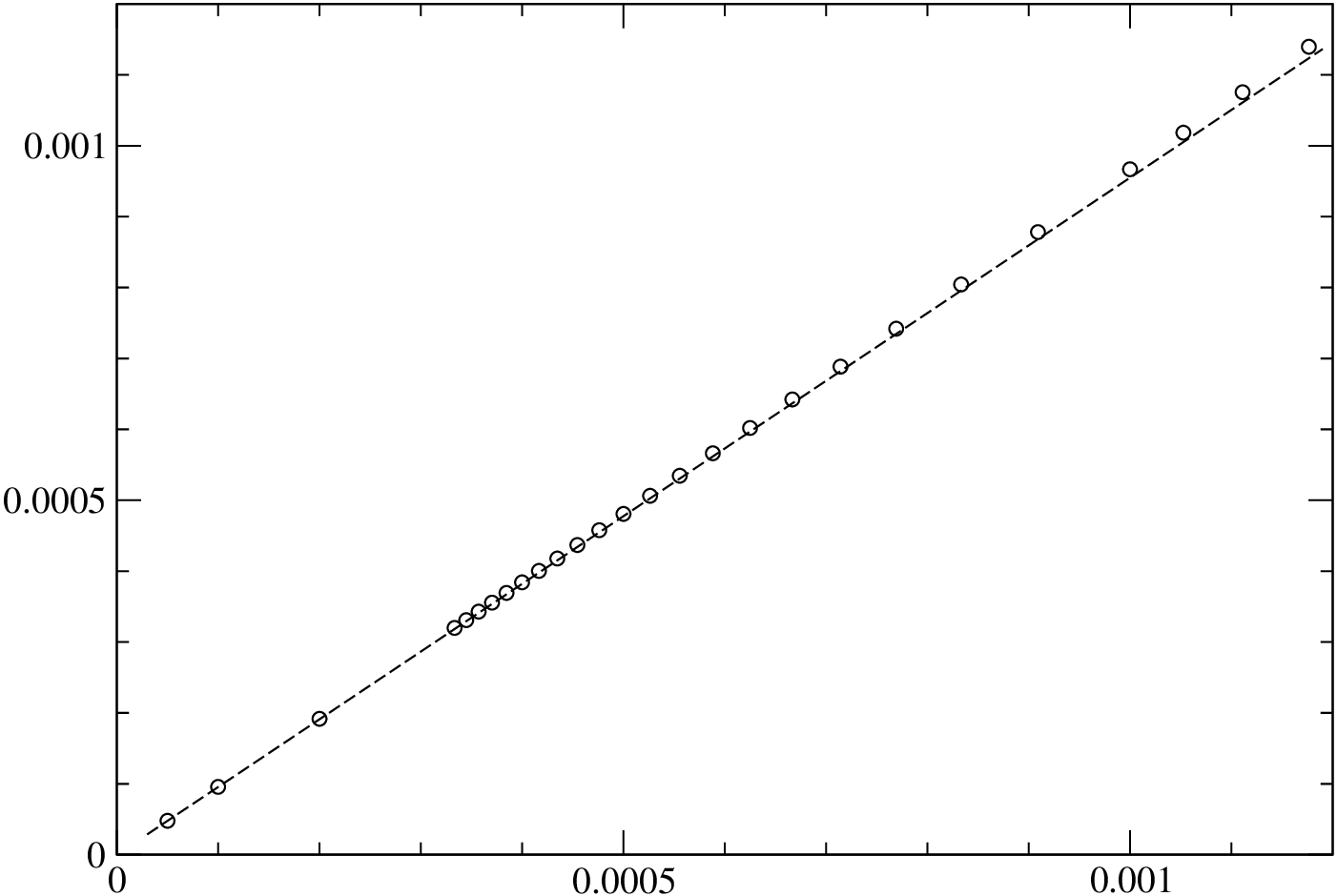}}}
\caption{The dependence of the $\rho_x$-component of real off-critical zeros
$(\rho_y=0)$ on the anisotropy parameter $\Delta$, calculated
numerically by using equation (\ref{exacteq}).
As explained in the text, $\rho_x=\frac{1}{2}$ at
$\Delta_c^*\approx 0.141733239663887$ and there are two conjugate solutions
$\rho_x$ and $1-\rho_x$ for $0<\Delta<\Delta_c^*$.
As $\Delta\to 0$, the two values of $\rho_x$ tend to the borders $0$
(the down branch) and $1$ (the up branch) of the critical strip.
The inset concerns the down branch and the region of small $\Delta$ where
numerical data for $\rho_x$ versus $\Delta$ (open circles) satisfy
the asymptotic relation $\rho_x\sim (3/\pi) \Delta$ (dashed line).}
\label{realoff}
\end{figure}

\textbf{Numerical observations.} The imaginary part of this point $\rho_y^*=0$ as well as its associated parameter $\Delta_c^*$ fulfill Eqs. (\ref{edge1}) and (\ref{edge22}) provided that
\begin{equation} \label{deltastar}
-2 + \int_0^1\left[
\theta_3\left({\rm e}^{-\pi t \Delta_c^*}\right)
\theta_3\left({\rm e}^{-\pi t/\Delta_c^*}\right) -\frac{1}{t} \right] \frac{{\rm d}t}{\sqrt{t}} 
= 0 .
\end{equation}
This equation obviously corresponds to $Z(1/2,\Delta_c^*)=0$ (see
\eqref{def:ZsDelta} and \eqref{final}).
Using the work of Montgomery \cite{Mont} (recently recovered by Faulhuber and
Steinerberger in \cite{FaulhuberSteinerberger17}),
we can show the following result.
\begin{lemma}
Equation \eqref{deltastar} admits a unique solution $\Delta_c^*$ on $(0,1]$.
\end{lemma}
\begin{proof}
It has been shown in \cite{FaulhuberSteinerberger17,Mont} that, for all $t>0$,
the function $\Delta\mapsto \theta_3\left({\rm e}^{-\pi t \Delta}\right)
\theta_3\left({\rm e}^{-\pi t/\Delta}\right)$ is strictly decreasing on $(0,1)$.
Furthermore, for $\Delta=1$ we have
$$
-2 + \int_0^1\left[ \theta_3\left(e^{-\pi t}\right)^2 -\frac{1}{t} \right]
\frac{{\rm d}t}{\sqrt{t}} \approx -1.9501325<0.
$$
To study the opposite $\Delta\to 0$ limit, we recall that
$Z(s,\Delta)$ is related to the Epstein zeta function
$\zeta^{(2)}(s,\Delta)$ via equation (\ref{def:ZsDelta}) and
the small-$\Delta$ behavior of $\zeta^{(2)}(s,\Delta)$ is given by
the asymptotic relation
(\ref{deltasmall2}) where the $p=2$ error term of the order
${\cal O}(\Delta)$ can be neglected in the limit $\Delta\to 0$.
Thus one arrives at the asymptotic relation
$$
Z(s,\Delta) \mathop{\sim}_{\Delta\to 0} \frac{1}{\Delta^s}
\frac{1}{\pi^{\frac{1}{2}-s}} \Gamma\left( \frac{1}{2}-s\right)
\zeta(1-2s) + \frac{1}{\Delta^{1-s}}
\pi^{\frac{1}{2}-s} \Gamma\left( s-\frac{1}{2}\right)
\zeta(2s-1)
$$
which exhibits the required duality symmetry (\ref{duality:ZsDelta}).
Consequently,
$$
\lim_{s\to\frac{1}{2}} Z(s,\Delta) =
\frac{\gamma-\log(4\pi\Delta)}{\sqrt{\Delta}} ,
$$
where $\gamma=0.5772156\ldots$ is the Euler-Mascheroni constant.
Thus,
$$
Z\left(\frac{1}{2},\Delta\right) \mathop{\sim}_{\Delta\to 0}
- \frac{\log\Delta}{\sqrt{\Delta}} \to + \infty .
$$
It follows that $\Delta_c^*$ is unique and the proof is complete.
\end{proof}
As is seen in Table \ref{edgetable}, the numerical solution of this
equation is $\Delta_c^*\approx 0.141733239663887$.
According to our numerical observations:
\begin{enumerate}
\item it turns out that the curve of off-critical zeros going out of
the critical edge zero 1 stays on the $\Delta$-axis, i.e. $\rho_y=0$.
\item As concerns the coupled equations for off-critical zeros (\ref{offzeros1})
and (\ref{offzeros2}), the second one is automatically fulfilled for
$\rho_y=0$ while the first one implies
\begin{equation} \label{exacteq}
- \left( \frac{1}{1+2\delta\rho_x} + \frac{1}{1-2\delta\rho_x} \right)
+ \int_0^1 \cosh\left( \delta\rho_x\log t\right)
\left[ \theta_3\left({\rm e}^{-\pi t \Delta}\right)
\theta_3\left({\rm e}^{-\pi t/\Delta}\right) -\frac{1}{t} \right]
\frac{{\rm d}t}{\sqrt{t}} = 0 .
\end{equation}
This equation has two conjugate \emph{real} solutions for $\delta\rho_x$
and $-\delta\rho_x$ (or, equivalently, $\rho_x$ and $1-\rho_x$),
only if $0\le \Delta<\Delta_c^*$ as observed in Figure \ref{realoff}.
\item The value of $\rho_x$ is $\frac{1}{2}$ at the critical edge point 1
corresponding to $\Delta=\Delta_c^*$.
As observed again in  Figure \ref{realoff}, decreasing the value
of $\Delta<\Delta_c^*$, the two values of $\rho_x$ split and tend
to the borders $0$ (down branch) and $1$ (up branch) of the critical strip
in the limit $\Delta\to 0$.
\end{enumerate}

\textbf{Heuristic derivation of $\Delta_c^*$ from small $\Delta$ approximation.}
It is useful to test how Eq. (\ref{zz}), which is accurate but certainly
only approximate for complex zeros (see Table \ref{approxtable}), works
in the present case of real off-critical zeros.
Writing $\rho=\frac{1}{2}+\delta\rho_x$ in (\ref{zz}), one gets
\begin{equation} \label{approx}
\left( \frac{\Delta}{\pi} \right)^{2\delta\rho_x} = -
\frac{\Gamma\left( -\delta\rho_x \right) \zeta\left( -2\delta\rho_x \right)}{
\Gamma\left( \delta\rho_x \right) \zeta\left( 2\delta\rho_x \right)} .      
\end{equation}
The expansion of both sides of this equation to the first order in small
$\delta\rho_x$ must be consistent at $\Delta=\Delta_c^*$:
as $\delta\rho_x\to 0$ we have
\begin{equation}
\left( \frac{\Delta^*}{\pi} \right)^{2\delta\rho_x}
= 1 + 2 \log\left( \frac{\Delta_c^*}{\pi} \right)
\delta\rho_x + \mathcal{O}\left( \delta\rho_x^2 \right) ,
\end{equation}
\begin{equation}
- \frac{\Gamma\left( -\delta\rho_x \right) \zeta\left( -2\delta\rho_x \right)}{
\Gamma\left( \delta\rho_x \right) \zeta\left( 2\delta\rho_x \right)}
= 1 + 2 \left[ \gamma-2\log(2\pi) \right]
\delta\rho_x + \mathcal{O}\left( \delta\rho_x^2 \right) .     
\end{equation}
Consequently, it must hold that
\begin{equation} \label{deltastaranalyt}
\Delta_c^* = \frac{e^\gamma}{4\pi} .  
\end{equation}  
We checked that $\Delta_c^*$ evaluated by using this analytic relation
coincides with the previous numerical estimate
$\Delta_c^*\approx0.141733239663887$ obtained from the exact
Eq. (\ref{deltastar}) by at least 22 decimal digits;
to go further a computer facility more powerful than the one at our
disposal is needed.
This indicates that the result (\ref{deltastaranalyt}) might be exact
which is difficult to prove directly by using (\ref{deltastar}). We have stated the corresponding open problem as Conjecture \ref{conjecture}.

\medskip

\textbf{Approximation in the $\rho_x,\Delta\to 0$ regime and comparison with numerical data.} As concerns the accuracy of the real off-critical zeros implied by
Eq. (\ref{approx}), for various values of $\Delta$ from the interval
$(0,\Delta_c^*)$ they coincide with the ones obtained from the exact Eq.
(\ref{exacteq}) up to 27 decimal digits, which supports the hypothesis
that the real off-critical zeros generated from Eq. (\ref{approx}) are exact.

For the down branch in Figure \ref{realoff}, to obtain the asymptotic
tendency of $\rho_x$ to $0$ as $\Delta\to 0$ by using the relation
(\ref{approx}), one writes $\delta\rho_x=\rho_x-\frac{1}{2}$ and expands
the right-hand side in small $\rho_x$:
\begin{equation}
- \frac{\Gamma\left( \frac{1}{2} - \rho_x\right) \zeta\left(1-2\rho_x\right)}{
\Gamma\left( \rho_x-\frac{1}{2}\right) \zeta\left(2\rho_x-1\right)}
=\frac{3}{\rho_x} +o(1),\quad \textnormal{as $\rho_x\to 0$}.
\end{equation}
Consequently, again as $\Delta\to 0$
\begin{equation}
\rho_x = \frac{3}{\pi} \Delta + o(1).
\end{equation}
As is seen in the inset of Figure \ref{realoff}, numerical data (open circles)
agree well with this analytic prediction (dashed line).

\begin{remark}
With regard to the symmetry $\Delta\to 1/\Delta$ of basic equations
for zeros of $\zeta^{(2)}(s,\Delta)$, there exists a pair of continuous
curves of real off-critical zeros also for each
$\Delta>1/\Delta_c^*\approx7.055507955448192$. 
\end{remark}

\section*{Acknowledgment}
Ladislav \v{S}amaj \and Igor Trav\v{e}nec acknowledge the support received from VEGA Grant No. 2/0092/21.

\end{document}